\documentclass[leqno,11pt]{amsart}
\usepackage{amsmath,amscd,amsthm}
\usepackage{epsfig,graphics,ifsym}
\usepackage{amssymb,latexsym}
\usepackage{mathrsfs,amsbsy}
\usepackage{hyperref}
\usepackage[all,cmtip]{xy}

\newtheorem{thm}{\bf Theorem}[section]

\newtheorem{prop}[thm]{\bf Proposition}
\newtheorem{cor}[thm]{\bf Corollary}
\newtheorem{lem}[thm]{\bf Lemma}
\newtheorem{rem}[thm]{\bf Remark}
\newtheorem{ex}[thm]{\bf Example}

\newtheorem{nono-theorem}{Theorem}[]
\newtheorem*{thm*}{Theorem}

\newcommand{\A}{\mathscr{A}}

\newcommand{\B}{\mathbf{B}}

\newcommand{\pf}{\noindent{\bfseries Proof. }}
\newcommand{\ov}{\overline}

\newcommand{\E}{\Lambda(V)}

\newcommand{\gl}{\mathfrak{gl}}
\newcommand{\Z}{\mathbb{Z}}

\newcommand{\te}{\widetilde{e}}
\newcommand{\tf}{\widetilde{f}}
\newcommand{\g}{\mathfrak{g}}
\newcommand{\td}{\widetilde}

\newcommand{\mf}{\mathfrak}

\newcommand{\psp}{\psi}
\newcommand{\psm}{\psi^\ast}

\newcommand{\om}{\omega}

\numberwithin{equation}{section}

\begin{document}
\title[ ]
{$q$-deformed Clifford algebra and level zero fundamental representations of quantum affine algebras}
\author{JAE-HOON KWON}
\address{ Department of Mathematics \\ Sungkyunkwan University \\ Suwon,  Republic of Korea}
\email{jaehoonkw@skku.edu}

\thanks{This work was  supported by Basic Science Research Program through the National Research Foundation of Korea (NRF)
funded by the Ministry of  Education, Science and Technology (No. 2012-0002607).}

\begin{abstract}
We give a realization of the level zero fundamental representations $W(\varpi_k)$ of the quantum affine algebra $U_q'(\mf{g})$, when $\mf{g}$ has a maximal parabolic subalgebra  of type $C_n$. 
We define a semisimple $U'_q({\mf g})$-module structure on $\E^{\otimes 2}$ in terms of $q$-deformed Clifford generators, where  $\E$ is the exterior algebra generated by the dual natural representation  $V$ of   $U_q(\mf{sl}_{n})$. We  show that each $W(\varpi_k)$ appears    in $\E^{\otimes 2}$  (not necessarily multiplicity-free). As a byproduct, we obtain a simple description of the crystal of $W(\varpi_k)$ in terms of $n\times 2$ binary matrices and their $(\mf{sl}_n,\mf{sl}_2)$-bicrystal structure.
\end{abstract}

\maketitle
\setcounter{tocdepth}{1}

\section{Introduction}
Let ${\mf g}$ be an affine Kac-Moody algebra and let $U'_q(\mf{g})$ be the associated quantum affine algebra without derivation. For a level zero fundamental weight  $\varpi_k$, Kashiwara introduced a finite-dimensional irreducible $U'_q(\mf{g})$-module $W(\varpi_k)$, which is called a level zero fundamental representation. It is obtained from a level zero extremal weight module $V(\varpi_k)$ by specializing its $U'_q({\mf g})$-linear automorphism $z_k$ as $1$, and it has a crystal base and global crystal base \cite{Kas02}. By the works of Chari and Pressley \cite{CP94,CP98},  any finite-dimensional irreducible $U'_q({\mf g})$-module is isomorphic to a subquotient of $W(\varpi_{k_1})_{a_1}\otimes \cdots \otimes W(\varpi_{k_r})_{a_r}$ for some $(k_1,a_1),\ldots,(k_r,a_r)$, where $W(\varpi_k)_a$ is obtained  by specializing $z_k$ as $a$. We also refer the reader to \cite{DO,JMO,KMN2,Ko} for previously known constructions of $W(\varpi_k)$ of various types.

The aim of this article is to introduce a  realization of $W(\varpi_k)$ and its crystal for a special class of affine Kac-Moody algebras ${\mf g}$, which  has a maximal parabolic subalgebra of type $C_n$, that is, ${\mf g}=C_{n}^{(1)}$, $A_{2n}^{(2)}$, $A_{2n}^{(2)\dagger}$, and $A_{2n-1}^{(2)}$. Instead of using $q$-wedge relations for classical Lie algebras of type $B$, $C$, and $D$, which were derived via $R$-matrix by Jing, Misra, and Okado \cite{JMO}, we construct a semisimple $U'_q({\mf g})$-module using a homomorphic image of $U'_q({\mf g})$ in a $q$-deformed Clifford algebra, which has a simple description of crystal structure and contains  all $W(\varpi_k)$ as its irreducible factor.

\begin{figure}
\includegraphics[width=7cm, height=20cm]{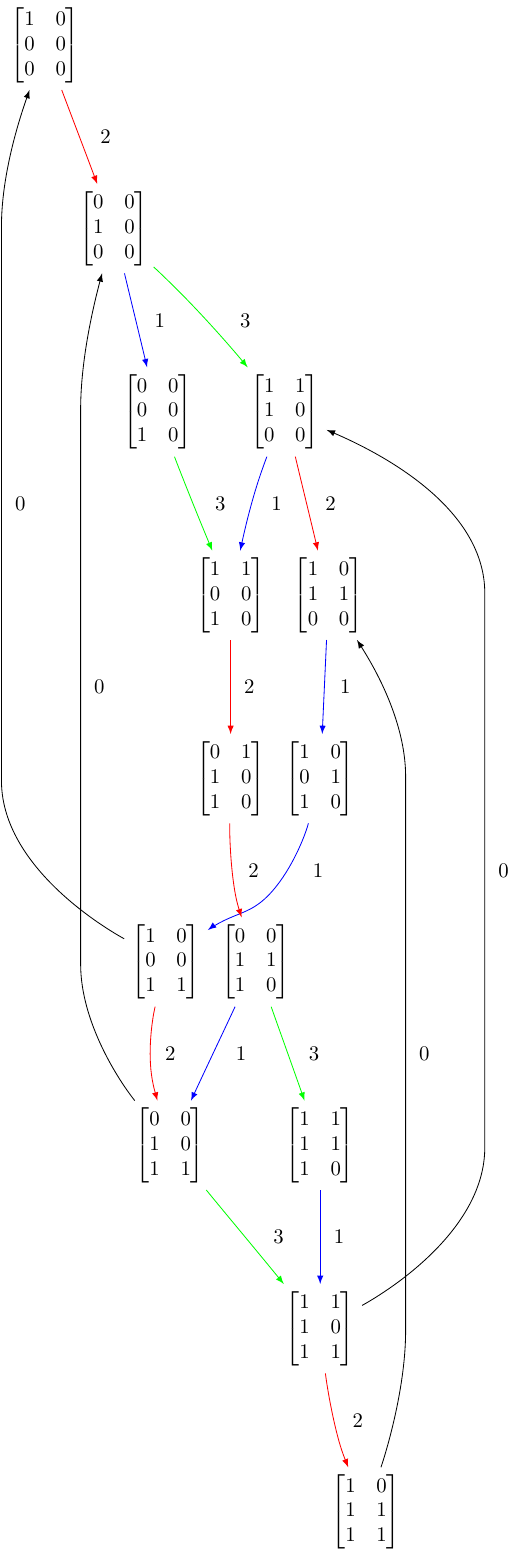}
\caption{ $B(W(\varpi_2))=C({\bf v}_{2,1})$ of type $C_3^{(1)}$. (This graph was implemented by SAGE.)} \label{Graph A}
\end{figure}

More precisely, we consider an exterior algebra  $\E$  generated by the dual natural representation  $V$ of   $U_q(\mf{sl}_{n})\subset U'_q({\mf g})$. Based on the action of $q$-deformed Clifford algebra on $\E$ due to Hayashi \cite{Ha}, we extend the $U_q(\mf{sl}_n)$-action on $\E^{\otimes 2}$ to that of $U'_q({\mf g})$, and show that it is a semisimple $U'_q({\mf g})$-module with a polarizable crystal base. 
The crystal of $\E^{\otimes 2}$ can be identified with the set of $n\times 2$ binary matrices, whose $U'_q({\mf g})$-crystal structure has a very simple description (see Figures \ref{Graph A}-\ref{Graph  C}).  Using its decomposition into connected components, we show that an irreducible summand in $\E^{\otimes 2}$ is generated by an extremal weight vector and then isomorphic to $W(\varpi_k)$ or $W(\varpi_k)_{-1}$ for some $k$ (Theorem \ref{main theorem}). We also obtain an explicit decomposition of $\E^{\otimes 2}$, where each $W(\varpi_k)$ appears at least once but not necessarily multiplicity-free (Corollary \ref{decomposition of E^2-2}). 

Moreover, from the $q$-deformed skew Howe duality on $\Lambda(\mathbb{C}^n\otimes \mathbb{C}^2)$ (cf. \cite{CKL,K13-3} and \cite{DK,La} for its crystal version), we observe that there are additional $U_q(\mf{sl}_2)$-crystal operators $\td{E}$ and $\td{F}$ acting on the crystal of $\E^{\otimes 2}$, which commute with those of $U_q(\mf{sl}_n)\subset U'_q({\mf g})$. This together with the author's previous work on classical crystals of type $B$ and $C$ \cite{K13-2} enables us to characterize the crystal of $W(\varpi_k)$ explicitly in terms of binary matrices and their statistics coming from an $\mf{sl}_2$-string with respect to $\td{E}$ and $\td{F}$ (Theorem \ref{main theorem-2}). The $(\mf{sl}_n,\mf{sl}_2)$-bicrystal structure also plays a crucial role in the decomposition of $\E^{\otimes 2}$.\vskip 2mm

{\bf Acknowledgment} The author would like to thank Myungho Kim for helpful discussion and the referee for careful reading.

\section{Background}
Let us briefly recall some necessary background on quantum affine algebras and crystal bases (see \cite{Kas02} for more details and references therein). 

\subsection{Notations} 
Let $A=(a_{ij})_{i,j\in I}$ be a generalized Cartan matrix
of affine type with an index set $I=\{\,0,1,\ldots,n\,\}$ and let $\mf{g}$ denote the associated affine Kac-Moody algebra with the Cartan subalgebra $\mf{h}$ \cite[\textsection 4.8]{K}.
Let $\{\,\alpha_i\,|\,i\in I\,\}\subset \mf{h}^*$ and $\{\,h_i\,|\,i\in I\,\}\subset \mf{h}$ be the sets of simple roots and simple coroots of $\mf{g}$, respectively, with $\langle h_i,\alpha_j\rangle =a_{ij}$ for $i,j\in I$.  We assume that $\{\,\alpha_i\,|\,i\in I\,\}$ and $\{\,h_i\,|\,i\in I\,\}$ are linearly independent. 
For $r\in \{0,n\}$, put $I_r=I\setminus\{r\}$ and let $\mf{g}_r$ be the subalgebra of $\mf{g}$ associated to $(a_{ij})_{i,j\in I_r}$. 

Let $c=\sum_{i\in I}a^\vee_ih_i$ be the canonical central element and let $\delta=\sum_{i\in I}a_i\alpha_i$ be the generator of the null roots. Let $\Lambda_i$ be the fundamental weight such that $\langle h_j,\Lambda_i \rangle=\delta_{ij}$ for $i,j\in I$.
We take a weight lattice $P$ such that $\alpha_i,\Lambda_i\in P$ and $h_i\in P^*={\rm Hom}(P,\Z)$. 
Let $(\ , \ )$ be a non-degenerate symmetric bilinear form on $\mf{h}^*$ satisfying $\langle h_i,\lambda \rangle =2(\alpha_i,\lambda)/(\alpha_i,\alpha_i)$ for $i\in I$ and $\lambda\in \mf{h}^*$, and normalized by $(\delta,\lambda)=\langle c,\lambda \rangle$ for $\lambda\in P$. 
Note that $(\alpha_i,\alpha_j)=a^\vee_ia_i^{-1}a_{ij}$ for $i,j\in I$. 

Let $ \mf{h}^*_{\rm cl}= \mf{h}^*/\mathbb{Q}\delta$ with the canonical projection ${\rm cl} : \mf{h}^* \longrightarrow \mf{h}^*_{\rm cl}$. Let $\mf{h}^{*0}=\{\,\lambda\in {\mf h}^*\,|\,\langle c,\lambda\rangle=0\,\}$ and $\mf{h}^{*0}_{\rm cl}={\rm cl}(\mf{h}^{*0})$. 
Let $P_{\rm cl}={\rm cl}(P)$, $P^0=\mf{h}^{*0}\cap P$, and $P^0_{\rm cl}={\rm cl}(P^0)$.

For $i\in I$, let $s_i$ be the simple reflection in $GL(\mf{h}^*)$ given by $s_i(\lambda)=\lambda-\langle h_i,\lambda\rangle \alpha_i$ for $\lambda\in \mf{h}^*$. Let $W$ be the Weyl group of $\mf{g}$ generated by $s_i$ for $i\in I$. Note that $W$ naturally induces an action on $\mf{h}^{*0}_{\rm cl}$, whose image we denote by $W_{\rm cl}$, and $W_{\rm cl}$ is generated by $s_i$ for $i\in I_0$.

\subsection{Quantum affine algebra and crystal base}
Let $d$ be the smallest positive integer such that $(\alpha_i,\alpha_i)/2 \in \frac{1}{d}\Z$ for $i\in I$. 
Let $q$ be an indeterminate and put $q_s=q^{1/d}$. Let $K=\mathbb{Q}(q_s)$. The quantum affine algebra $U_q(\mf{g})$ is the associative $K$-algebra with $1$ generated by $e_i$, $f_i$, and $q^h$ for $i\in I$ and $h\in \frac{1}{d}P^*$ subject to the relations:
{\allowdisplaybreaks \begin{gather*}
q^0=1, \quad q^{h +h'}=q^{h}q^{h'},\label{Rel-1} \\ q^h e_i=q^{\langle h,\alpha_i\rangle}
e_i q^h, \quad q^h f_i=q^{-\langle h,\alpha_i\rangle} f_i q^h, \label{Rel-2}\\ 
 e_i f_j-f_j e_i =\delta_{ij}\frac{t_i-t_i^{-1}} {q_i-q^{-1}_i},\label{Rel-3} \\ 
 \sum_{k=0}^{1-a_{ij}}(-1)^k e_i^{(k)} e_j e_i^{(1-a_{ij}-k)}=\sum_{k=0}^{1-a_{ij}}(-1)^k f_i^{(k)} f_j f_i^{(1-a_{ij}-k)}   =0 \ \ \ \ \ \ (i\neq j),\label{Rel-4}
\end{gather*}}
\hskip -2mm where $q_i = q^{(\alpha_i,\alpha_i)/2}$,  $t_i = q^{(\alpha_i,\alpha_i)h_i/2}$, and 
\begin{equation*}
[k]_i =\frac{q_i^k-q_i^{-k}}{q_i-q_i^{-1}},\quad
[k]_i!=\prod_{s=1}^{k}[s]_i, \quad
e_i^{(k)}=\frac{1}{[k]_i!}e_i^k, \quad f_i^{(k)}=\frac{1}{[k]_i!}f_i^k,
\end{equation*}
for $i\in I$ and $k\geq 0$.
Recall that $U_q(\mf{g})$ has a comultiplication $\Delta$ given by
\begin{equation*}\label{Delta_+}
\begin{split}
\Delta(q^h)&=q^h\otimes q^h, \\ \Delta(e_i)&=e_i\otimes t_i^{-1}+ 1\otimes e_i,\\ \Delta(f_i)&=f_i\otimes 1 + t_i\otimes f_i, \\  
\end{split}
\end{equation*}
for $i\in I$ and $h\in \frac{1}{d}P^*$.

We denote by $U_q'(\mf{g})$ the subalgebra of $U_q(\mf{g})$ generated by $e_i$, $f_i$, and $q^h$ for $i\in I$ and $h\in \frac{1}{d}(P_{\rm cl})^*$. Let $z$ be an indeterminate.
For a $U'_q(\mf{g})$-module $M$ with weight space decomposition $M=\bigoplus_{\lambda\in P_{\rm cl}}M_\lambda$, let $M_{\rm aff}=K[z,z^{-1}]\otimes M$ be a $U_q(\mf{g})$-module, where $e_i$ and $f_i\in U_q'(\mf{g})$ act by $z^{\delta_{0i}}\otimes e_i$ and $z^{-\delta_{0i}}\otimes f_i$, respectively  for $i\in I$, and ${\rm wt}(z^k\otimes m)={\rm wt}(m)+k\delta$ for $m\in M$ and $k\in\Z$.  Here ${\rm wt}$ denotes the weight function.
For $a\in K$, we define a $U'_q(\mf{g})$-module $M_a=M_{\rm aff}/(z-a)M_{\rm aff}$.

Let $M$ be an integrable module over $U_q(\mf{g})$ or $U'_q(\mf{g})$ having weight space decomposition $M=\bigoplus_{\lambda}M_\lambda$ with ${\rm dim}M_\lambda<\infty$ for $\lambda\in P$ or $P_{\rm cl}$. For $u\in M_{\lambda}$ and $i\in I$, we have
$u=\sum_{r\geq 0, -\langle h_i,\lambda \rangle}f_i^{(r)}u_r$, 
where $e_iu_r=0$ for all $r\ge 0$. 
We define $\te_i$ and $\tf_i$  by
$\te_i u=\sum_{r \ge 1} f_i^{(r-1)}u_r$ and     
$\tf_i u=\sum_{r \ge 0} f_i^{(r+1)}u_r$.
Let $\mathbb{A}$ denote the subring of $K$ consisting
of all rational functions which are regular at $q_s=0$. A pair $(L,B)$ is called a crystal base of $M$ if  
\begin{itemize}
\item[(1)] $L$ is an $\mathbb{A}$-lattice of $M$, where $L=\bigoplus_{\lambda}L_{\lambda}$ with $L_{\lambda}=L\cap M_{\lambda}$,
\item[(2)]  $\te_i L\subset L$ and $\tf_i L\subset L$ for $i\in I$,
\item[(3)] $B$ is a $\mathbb{Q}$-basis of $L/q_sL$, where $B=\bigsqcup_{\lambda}B_{\lambda}$ with
$B_{\lambda}=B\cap(L/q_sL)_{\lambda}$,

\item[(4)] $\te_iB \subset B\sqcup \{0\}$,
$\tf_i B\subset B\sqcup \{0\}$ for $i\in I$,

\item[(5)] for $b,b'\in B$ and $i\in I$,  
$\tf_i b = b'$ if and only if $b=\te_i b'$.
\end{itemize}

Following \cite{BKK} (cf. \cite{KMN2}), we say that a symmetric bilinear form $(\ , \ )$ on $M$ is a polarization if 
\begin{equation}\label{polarization}
(xu, v)=(u,\eta(x)v),
\end{equation}
for $x\in U_q(\mf{g})$ or $U'_q(\mf{g})$ and $u, v\in M$, where $\eta$ is the anti-automorphism given by
\begin{equation*}
\eta(q^h)=q^h, \ \ \eta(e_i)=q_i^{-1}t_i^{-1}f_i , \ \ \eta(f_i)=q_i^{-1}t_ie_i  \ \ \ (i\in I),
\end{equation*}
and say that a crystal base $(L,B)$ of $M$ is polarizable if $(L,L)\subset \mathbb{A}$ with respect to a polarization on $M$ and  $B$ is orthonormal (up to scalar multiplication by $\pm 1$) with respect to the induced $\mathbb{Q}$-bilinear form $(\ ,\ )_0$ on $L/q_sL$.  If $(\ ,\ )_{M_i}$ is a polarization of $M_i$ $(i=1,2)$, then $M_1\otimes M_2$ has a polarization given by $(u_1\otimes u_2,v_1\otimes v_2)_{M_1\otimes M_2}=(u_1,v_1)_{M_1}(u_2,v_2)_{M_2}$ for $u_i, v_i \in M_i$.  If $(L_i,B_i)$ is a polarizable crystal base of $M_i$, then $(L_1\otimes L_2, B_1\otimes B_2)$ is a polarizable crystal base of $M_1\otimes M_2$.

\begin{prop}[Theorem 2.12 in \cite{BKK}]\label{polarizable implies semisimple}
If $M$ has a polarizable crystal base, then $M$ is completely reducible.
\end{prop}

\subsection{Level zero fundamental representation}
For a regular crystal $B$, we define the action of $W$ as follows. For $i\in I$ and $b\in B$,
\begin{equation*}
{\texttt S}_{s_i}(b)=
\begin{cases}
\tf_i^{\, \langle h_i, {\rm wt}(b)\rangle }b, & \text{if $\langle h_i, {\rm wt}(b)\rangle \geq 0$},\\
\te_i^{\,-\langle h_i, {\rm wt}(b)\rangle }b, & \text{if $\langle h_i, {\rm wt}(b)\rangle \leq 0$}.\\
\end{cases}
\end{equation*}
For $w\in W$ with  a reduced expression $w=s_{i_1}\ldots s_{i_r}$, we let ${\texttt S}_w={\texttt S}_{s_{i_1}}\ldots {\texttt S}_{s_{i_r}}$.

Let $u_\lambda$ be a weight vector of an integrable module $M$ over $U_q(\mf{g})$ or  $U'_q(\mf{g})$ with  weight $\lambda$. Then $u_\lambda$ is called an extremal weight vector of extremal weight $\lambda$ if there exists $\{u_{w\lambda}\}_{w\in W}$ such that
\begin{itemize}
\item[(1)] $u_{w\lambda}=u_\lambda$ if $w$ is the identity,

\item[(2)] $e_i u_{w\lambda}=0$ and $f_i^{(\langle h_i,w\lambda\rangle)}u_{w\lambda}=u_{s_iw\lambda}$ if $\langle h_i,w\lambda\rangle\geq 0$,

\item[(3)] $f_i u_{w\lambda}=0$ and $e_i^{(-\langle h_i,w\lambda\rangle)}u_{w\lambda}=u_{s_iw\lambda}$ if $\langle h_i,w\lambda\rangle\leq 0$.
\end{itemize} 
If $u_{\lambda}$ is an extremal weight vector, we denote $u_{w\lambda}$ by $S_w u_\lambda$ for $w\in W$. 

For $\lambda\in P$, define $V(\lambda)$ to be a $U_q(\mf{g})$-module generated by a vector $u_\lambda$ of weight $\lambda$ subject to the relations such that $u_\lambda$ is an extremal weight vector. We call $V(\lambda)$ an extremal weight module with extremal weight $\lambda$. The notion of extremal weight module was introduced by Kashiwara and it was proved that $V(\lambda)$ has a crystal base and a global crystal base \cite{Kas94'}. Note that ${\texttt S}_w u_\lambda \equiv S_w u_\lambda \pmod{qL(\lambda)}$ for $w\in W$.

For $i\in I_0$, let
\begin{equation*}
\varpi_i =\Lambda_i- \langle c, \Lambda_i\rangle \Lambda_0=\Lambda_i- a_i^\vee \Lambda_0\in P^0,
\end{equation*}
which is called a level zero fundamental weight. We have 
$\{\,k\in\Z\,|\,\varpi_i+k\delta\in W\varpi_i\,\}=\Z d_i$, where $d_i=\max\{1,(\alpha_i,\alpha_i)/2\}$ except in the case $d_i=1$ for $\mf{g}=A_{2n}^{(2)}$ and $i=n$.
There exists a $U_q'(\mf{g})$-linear automorphism $z_i$ on $V(\varpi_i)$ of weight $d_i\delta$ sending $u_{\varpi_i}$ to $u_{\varpi_i+d_i\delta}$.  We define 
\begin{equation*}
W(\varpi_i) = V(\varpi_i)/(z_i-1)V(\varpi_i),
\end{equation*}
which is called a level zero fundamental representation of $U_q'(\mf{g})$ \cite{Kas02}. They play a crucial role, especially as building blocks of finite-dimensional $U'_q(\mf{g})$-modules. The following properties of $W(\varpi_i)$ are known, which is a part of \cite[Theorem 5.17]{Kas02}.

\begin{thm}\label{Fundamental Repn}\mbox{}
\begin{itemize}
\item[(1)] $W(\varpi_i)$ is a finite-dimensional irreducible integrable $U'_q(\mf{g})$-module.

\item[(2)] $W(\varpi_i)$ has a global crystal base with a simple crystal.

\item[(3)] $\dim W(\varpi_i)_\mu =1$ for $\mu\in W{\rm cl}(\varpi_i)$.

\item[(4)] The weight of an extremal weight vector of $W(\varpi_i)$ is in $W{\rm cl}(\varpi)$.

\item[(5)] The set of weights of $W(\varpi_i)$ is the intersection of ${\rm cl}(\varpi_i +\sum_{i\in I}\Z\alpha_i)$ and the convex hull of $W{\rm cl}(\varpi_i)$.

\item[(6)] Any finite-dimensional irreducible integrable $U'_q(\mf{g})$-module with ${\rm cl}(\varpi_i)$ as an extremal weight is isomorphic to $W(\varpi_i)_a$ for some $a\in K\setminus\{0\}$.
\end{itemize}
\end{thm}

\section{wedge spaces and level zero representations}
\subsection{Non-exceptional affine algebras} Throughout this paper, we assume that $\mf{g}$ is of type 
$B_n^{(1)}$, $C_n^{(1)}$, $D_n^{(1)}$, $A_{2n}^{(2)}$, $A_{2n-1}^{(2)}$, $D_{n+1}^{(2)}$ (called non-exceptional affine type together with $A_n^{(1)}$)  
following \cite{K} for the labeling of simple roots. Note that $\mf{g}_0 \cap  \mf{g}_n=A_{n-1}$ and $\mf{g}_r$ ($r=0,n$) is one of $B_n$, $C_n$, and $D_n$. 
Let us denote the type of $\mf{g}_r$ by a partition or a Young diagram $\diamondsuit$ as follows:
\begin{equation*}
\diamondsuit=
\begin{cases}
\ {\def\lr#1{\multicolumn{1}{|@{\hspace{.6ex}}c@{\hspace{.6ex}}|}{\raisebox{-.3ex}{$#1$}}}\raisebox{-.3ex}
{$\scalebox{0.45}{\begin{array}[b]{c}
\cline{1-1}
\lr{\ \ \ }\\ 
\cline{1-1}
\end{array}}$}}\ =(1) & \text{for $B_n$},\\
{\def\lr#1{\multicolumn{1}{|@{\hspace{.6ex}}c@{\hspace{.6ex}}|}{\raisebox{-.08ex}{$#1$}}}\raisebox{-.2ex}
{$\scalebox{0.45}{\begin{array}[b]{cc}
\cline{1-1}\cline{2-2}
\lr{\ \ \ }&\lr{\ \ \ }\\ 
\cline{1-1}\cline{2-2}
\end{array}}$}}=(2) & \text{for $C_n$}, \\
\ {\def\lr#1{\multicolumn{1}{|@{\hspace{.6ex}}c@{\hspace{.6ex}}|}{\raisebox{-.3ex}{$#1$}}}\raisebox{-1.2ex}
{$\scalebox{0.45}{\begin{array}[b]{c}
\cline{1-1}
\lr{\ \ \ }\\ 
\cline{1-1}
\lr{\ \ \ }\\
\cline{1-1}
\end{array}}$}}\ =(1,1) & \text{for $D_n$}.
\end{cases}
\end{equation*} \vskip 2mm

Since the type of $\mf{g}$ is completely determined by those of $\mf{g}_r$ ($r=0,n$), we may identify the type of $\mf{g}$ with a pair of partitions  $\pmb{\diamondsuit}=(\diamondsuit_0, \diamondsuit_n)$, where $\diamondsuit_r$ is the type of $\mf{g}_{n-r}$.  Since $\mf{g}_0 \cap  \mf{g}_n=A_{n-1}$ is fixed,  we may understand that $\diamondsuit_r$ is determined by $\alpha_r$ for $r\in\{0,n\}$  in the Dynkin diagram of $\mf{g}$. This convention will be useful when we realize $W(\varpi_i)$ and its crystal in later sections.
For the readers' convenience, we list the diagrams of $\mf{g}$ and the associated pair $\pmb{\diamondsuit}$.
\vskip 7mm

%B_n^{(1)}
{\small 
\begin{center}\setlength{\unitlength}{0.15in} \medskip
\ \, \scalebox{0.8}{\begin{picture}(28,4)
\put(0,2){\makebox(0,0)[c]{$B_n^{(1)}$}}
\put(6,0){\makebox(0,0)[c]{$\bigcirc$}}
\put(6,4){\makebox(0,0)[c]{$\bigcirc$}}
\put(8,2){\makebox(0,0)[c]{$\bigcirc$}}

\put(14,2){\makebox(0,0)[c]{$\bigcirc$}}
\put(16.2,2){\makebox(0,0)[c]{$\bigcirc$}}

\put(6.35,0.3){\line(1,1){1.35}} 
\put(6.35,3.7){\line(1,-1){1.35}}
\put(8.4,2){\line(1,0){1.55}} 
\put(12,2){\line(1,0){1.55}}
\put(11,1.95){\makebox(0,0)[c]{$\cdots$}}
\put(15.1,1.95){\makebox(0,0)[c]{$\Longrightarrow$}}

\put(6,5){\makebox(0,0)[c]{\tiny $0$}}
\put(6,-1){\makebox(0,0)[c]{\tiny $1$}}
\put(8,1){\makebox(0,0)[c]{\tiny $2$}}
\put(13.9,1){\makebox(0,0)[c]{\tiny ${n-1}$}}
\put(16.2,1){\makebox(0,0)[c]{\tiny ${n}$}}

\put(23,2){\makebox(0,0)[c]{$(\ {\def\lr#1{\multicolumn{1}{|@{\hspace{.6ex}}c@{\hspace{.6ex}}|}{\raisebox{-.3ex}{$#1$}}}\raisebox{-.8ex}
{$\scalebox{0.4}{\begin{array}[b]{c}
\cline{1-1}
\lr{\ \ \ }\\ 
\cline{1-1}
\lr{\ \ \ }\\
\cline{1-1}
\end{array}}$}}\ ,\ {\def\lr#1{\multicolumn{1}{|@{\hspace{.6ex}}c@{\hspace{.6ex}}|}{\raisebox{-.3ex}{$#1$}}}\raisebox{-.3ex}
{$\scalebox{0.4}{\begin{array}[b]{c}
\cline{1-1}
\lr{\ \ \ }\\ 
\cline{1-1}
\end{array}}$}} \ )$}}

\end{picture}}
\end{center}\vskip 2mm

%C_n^{(1)}

\begin{center}
\hskip -.7cm  \setlength{\unitlength}{0.15in}
\scalebox{0.8}{\begin{picture}(25,4)
\put(0,2){\makebox(0,0)[c]{$C_{n}^{(1)}$}}

\put(5.6,2){\makebox(0,0)[c]{$\bigcirc$}}
\put(8,2){\makebox(0,0)[c]{$\bigcirc$}}
\put(14.25,2){\makebox(0,0)[c]{$\bigcirc$}}
\put(16.5,2){\makebox(0,0)[c]{$\bigcirc$}}

\put(8.35,2){\line(1,0){1.5}}
\put(12.28,2){\line(1,0){1.45}}
\put(15.4,2){\makebox(0,0)[c]{$\Longleftarrow$}}

\put(6.8,1.97){\makebox(0,0)[c]{$\Longrightarrow$}}
\put(11,1.95){\makebox(0,0)[c]{$\cdots$}}

\put(5.6,1){\makebox(0,0)[c]{\tiny $0$}}
\put(8,1){\makebox(0,0)[c]{\tiny $1$}}
\put(14.2,1){\makebox(0,0)[c]{\tiny ${n-1}$}}
\put(16.5,1){\makebox(0,0)[c]{\tiny ${n}$}}

\put(22.9,2){\makebox(0,0)[c]{$(\ {\def\lr#1{\multicolumn{1}{|@{\hspace{.6ex}}c@{\hspace{.6ex}}|}{\raisebox{-.08ex}{$#1$}}}\raisebox{-.3ex}
{$\scalebox{0.4}{\begin{array}[b]{cc}
\cline{1-1}\cline{2-2}
\lr{\ \ \ }&\lr{\ \ \ }\\ 
\cline{1-1}\cline{2-2}
\end{array}}$}} \ ,\ {\def\lr#1{\multicolumn{1}{|@{\hspace{.6ex}}c@{\hspace{.6ex}}|}{\raisebox{-.08ex}{$#1$}}}\raisebox{-.3ex}
{$\scalebox{0.4}{\begin{array}[b]{cc}
\cline{1-1}\cline{2-2}
\lr{\ \ \ }&\lr{\ \ \ }\\ 
\cline{1-1}\cline{2-2}
\end{array}}$}}  \ )$}}
\end{picture}}
\end{center}\vskip 2mm

%D_n^{(1)}

\begin{center}\setlength{\unitlength}{0.15in}\medskip
\hskip -0.5cm \scalebox{0.8}{\begin{picture}(25,4)
\put(0,2){\makebox(0,0)[c]{$D_n^{(1)}$  }}
\put(6,0){\makebox(0,0)[c]{$\bigcirc$}}
\put(6,4){\makebox(0,0)[c]{$\bigcirc$}}
\put(8,2){\makebox(0,0)[c]{$\bigcirc$}}

\put(14,2){\makebox(0,0)[c]{$\bigcirc$}}
\put(16,0){\makebox(0,0)[c]{$\bigcirc$}}
\put(16,4){\makebox(0,0)[c]{$\bigcirc$}}

\put(6.35,0.3){\line(1,1){1.35}} \put(6.35,3.7){\line(1,-1){1.35}}
\put(8.4,2){\line(1,0){1.55}} \put(12,2){\line(1,0){1.55}}

\put(11,1.95){\makebox(0,0)[c]{$\cdots$}}

\put(14.35,2.3){\line(1,1){1.35}} \put(14.35,1.65){\line(1,-1){1.35}}

\put(6,5){\makebox(0,0)[c]{\tiny $0$}}
\put(6,-1){\makebox(0,0)[c]{\tiny $1$}}
\put(8,1){\makebox(0,0)[c]{\tiny $2$}}
\put(13.7,1){\makebox(0,0)[c]{\tiny ${n-2}$}}
\put(16,5){\makebox(0,0)[c]{\tiny ${n-1}$}}
\put(16,-1){\makebox(0,0)[c]{\tiny ${n}$}}

\put(22.8,2){\makebox(0,0)[c]{$(\ {\def\lr#1{\multicolumn{1}{|@{\hspace{.6ex}}c@{\hspace{.6ex}}|}{\raisebox{-.3ex}{$#1$}}}\raisebox{-.8ex}
{$\scalebox{0.4}{\begin{array}[b]{c}
\cline{1-1}
\lr{\ \ \ }\\ 
\cline{1-1}
\lr{\ \ \ }\\
\cline{1-1}
\end{array}}$}}\ ,\ {\def\lr#1{\multicolumn{1}{|@{\hspace{.6ex}}c@{\hspace{.6ex}}|}{\raisebox{-.3ex}{$#1$}}}\raisebox{-.8ex}
{$\scalebox{0.4}{\begin{array}[b]{c}
\cline{1-1}
\lr{\ \ \ }\\ 
\cline{1-1}
\lr{\ \ \ }\\
\cline{1-1}
\end{array}}$}} \ )$}}

\end{picture}}
\end{center}\vskip 2mm

\begin{center}
\hskip -  .7cm  \setlength{\unitlength}{0.15in}
\scalebox{0.8}{\begin{picture}(25,4)
\put(0,2){\makebox(0,0)[c]{$A_{2n}^{(2)}$}}

\put(5.6,2){\makebox(0,0)[c]{$\bigcirc$}}
\put(8,2){\makebox(0,0)[c]{$\bigcirc$}}
\put(14.25,2){\makebox(0,0)[c]{$\bigcirc$}}
\put(16.5,2){\makebox(0,0)[c]{$\bigcirc$}}

\put(8.35,2){\line(1,0){1.5}}
\put(12.28,2){\line(1,0){1.45}}
\put(15.4,2){\makebox(0,0)[c]{$\Longleftarrow$}}

\put(6.8,1.97){\makebox(0,0)[c]{$\Longleftarrow$}}
\put(11,1.95){\makebox(0,0)[c]{$\cdots$}}

\put(5.6,1){\makebox(0,0)[c]{\tiny $0$}}
\put(8,1){\makebox(0,0)[c]{\tiny $1$}}
\put(14.2,1){\makebox(0,0)[c]{\tiny ${n-1}$}}
\put(16.5,1){\makebox(0,0)[c]{\tiny ${n}$}}

\put(23,2){\makebox(0,0)[c]{$(\ {\def\lr#1{\multicolumn{1}{|@{\hspace{.6ex}}c@{\hspace{.6ex}}|}{\raisebox{-.3ex}{$#1$}}}\raisebox{-.3ex}
{$\scalebox{0.4}{\begin{array}[b]{c}
\cline{1-1}
\lr{\ \ \ }\\ 
\cline{1-1}
\end{array}}$}}\ ,\ {\def\lr#1{\multicolumn{1}{|@{\hspace{.6ex}}c@{\hspace{.6ex}}|}{\raisebox{-.08ex}{$#1$}}}\raisebox{-.3ex}
{$\scalebox{0.4}{\begin{array}[b]{cc}
\cline{1-1}\cline{2-2}
\lr{\ \ \ }&\lr{\ \ \ }\\ 
\cline{1-1}\cline{2-2}
\end{array}}$}} \ )$}}

\end{picture}}
\end{center}\vskip 2mm

\begin{center}
\hskip -  .7cm  \setlength{\unitlength}{0.15in}
\scalebox{0.8}{\begin{picture}(25,4)
\put(0,2){\makebox(0,0)[c]{$A_{2n}^{(2)\dagger}$}}

\put(5.6,2){\makebox(0,0)[c]{$\bigcirc$}}
\put(8,2){\makebox(0,0)[c]{$\bigcirc$}}
\put(14.25,2){\makebox(0,0)[c]{$\bigcirc$}}
\put(16.5,2){\makebox(0,0)[c]{$\bigcirc$}}

\put(8.35,2){\line(1,0){1.5}}
\put(12.28,2){\line(1,0){1.45}}
\put(15.4,2){\makebox(0,0)[c]{$\Longrightarrow$}}

\put(6.8,1.97){\makebox(0,0)[c]{$\Longrightarrow$}}
\put(11,1.95){\makebox(0,0)[c]{$\cdots$}}

\put(5.6,1){\makebox(0,0)[c]{\tiny $0$}}
\put(8,1){\makebox(0,0)[c]{\tiny $1$}}
\put(14.2,1){\makebox(0,0)[c]{\tiny ${n-1}$}}
\put(16.5,1){\makebox(0,0)[c]{\tiny ${n}$}}

\put(23,2){\makebox(0,0)[c]{$(\ {\def\lr#1{\multicolumn{1}{|@{\hspace{.6ex}}c@{\hspace{.6ex}}|}{\raisebox{-.08ex}{$#1$}}}\raisebox{-.3ex}
{$\scalebox{0.4}{\begin{array}[b]{cc}
\cline{1-1}\cline{2-2}
\lr{\ \ \ }&\lr{\ \ \ }\\ 
\cline{1-1}\cline{2-2}
\end{array}}$}}\ ,\ {\def\lr#1{\multicolumn{1}{|@{\hspace{.6ex}}c@{\hspace{.6ex}}|}{\raisebox{-.3ex}{$#1$}}}\raisebox{-.3ex}
{$\scalebox{0.4}{\begin{array}[b]{c}
\cline{1-1}
\lr{\ \ \ }\\ 
\cline{1-1}
\end{array}}$}} \ )$}}

\end{picture}}
\end{center}\vskip 2mm

\begin{center}\setlength{\unitlength}{0.15in} \medskip
\hskip -.5cm \scalebox{0.8}{\begin{picture}(25,4)
\put(0,2){\makebox(0,0)[c]{$A_{2n-1}^{(2)}$}}
\put(6,0){\makebox(0,0)[c]{$\bigcirc$}}
\put(6,4){\makebox(0,0)[c]{$\bigcirc$}}
\put(8,2){\makebox(0,0)[c]{$\bigcirc$}}

\put(14,2){\makebox(0,0)[c]{$\bigcirc$}}
\put(16.2,2){\makebox(0,0)[c]{$\bigcirc$}}

\put(6.35,0.3){\line(1,1){1.35}} 
\put(6.35,3.7){\line(1,-1){1.35}}
\put(8.4,2){\line(1,0){1.55}} 
\put(12,2){\line(1,0){1.55}}
\put(11,1.95){\makebox(0,0)[c]{$\cdots$}}
\put(15.1,1.95){\makebox(0,0)[c]{$\Longleftarrow$}}

\put(6,5){\makebox(0,0)[c]{\tiny $0$}}
\put(6,-1){\makebox(0,0)[c]{\tiny $1$}}
\put(8,1){\makebox(0,0)[c]{\tiny $2$}}
\put(13.9,1){\makebox(0,0)[c]{\tiny ${n-1}$}}
\put(16.2,1){\makebox(0,0)[c]{\tiny ${n}$}}

\put(23,2){\makebox(0,0)[c]{$(\ {\def\lr#1{\multicolumn{1}{|@{\hspace{.6ex}}c@{\hspace{.6ex}}|}{\raisebox{-.3ex}{$#1$}}}\raisebox{-.8ex}
{$\scalebox{0.4}{\begin{array}[b]{c}
\cline{1-1}
\lr{\ \ \ }\\ 
\cline{1-1}
\lr{\ \ \ }\\
\cline{1-1}
\end{array}}$}}\ ,\ {\def\lr#1{\multicolumn{1}{|@{\hspace{.6ex}}c@{\hspace{.6ex}}|}{\raisebox{-.08ex}{$#1$}}}\raisebox{-.3ex}
{$\scalebox{0.4}{\begin{array}[b]{cc}
\cline{1-1}\cline{2-2}
\lr{\ \ \ }&\lr{\ \ \ }\\ 
\cline{1-1}\cline{2-2}
\end{array}}$}} \ )$}}

\end{picture}}
\end{center}\vskip 2mm

\begin{center}
\hskip -.7cm  \setlength{\unitlength}{0.15in}
\scalebox{0.8}{\begin{picture}(25,4)
\put(0,2){\makebox(0,0)[c]{$D_{n+1}^{(2)}$}}

\put(5.6,2){\makebox(0,0)[c]{$\bigcirc$}}
\put(8,2){\makebox(0,0)[c]{$\bigcirc$}}
\put(14.25,2){\makebox(0,0)[c]{$\bigcirc$}}
\put(16.5,2){\makebox(0,0)[c]{$\bigcirc$}}

\put(8.35,2){\line(1,0){1.5}}
\put(12.28,2){\line(1,0){1.45}}
\put(15.4,2){\makebox(0,0)[c]{$\Longrightarrow$}}

\put(6.8,1.97){\makebox(0,0)[c]{$\Longleftarrow$}}
\put(11,1.95){\makebox(0,0)[c]{$\cdots$}}

\put(5.6,1){\makebox(0,0)[c]{\tiny $0$}}
\put(8,1){\makebox(0,0)[c]{\tiny $1$}}

\put(14.3,1){\makebox(0,0)[c]{\tiny ${n-1}$}}
\put(16.5,1){\makebox(0,0)[c]{\tiny ${n}$}}

\put(23,2){\makebox(0,0)[c]{$(\ {\def\lr#1{\multicolumn{1}{|@{\hspace{.6ex}}c@{\hspace{.6ex}}|}{\raisebox{-.3ex}{$#1$}}}\raisebox{-.3ex}
{$\scalebox{0.4}{\begin{array}[b]{c}
\cline{1-1}
\lr{\ \ \ }\\ 
\cline{1-1}
\end{array}}$}}\ ,\ {\def\lr#1{\multicolumn{1}{|@{\hspace{.6ex}}c@{\hspace{.6ex}}|}{\raisebox{-.3ex}{$#1$}}}\raisebox{-.3ex}
{$\scalebox{0.4}{\begin{array}[b]{c}
\cline{1-1}
\lr{\ \ \ }\\ 
\cline{1-1}
\end{array}}$}} \ )$}}
\end{picture}}
\end{center}}
Here $A_{2n}^{(2)\dagger}$ is a different labeling of simple roots of $A_{2n}^{(2)}$, and in this case, we have
\begin{equation*}
\varpi_n=2\Lambda_n-\Lambda_0, \ \ \varpi_{i}=\Lambda_i-\Lambda_0 \ \ (i=1,\ldots,n-1).
\end{equation*}

\subsection{$q$-deformed Clifford algebra}
Let $[\ov{n}]=\{\ov{n}<\ldots<\ov{1}\}$ be a linearly ordered set. 
Consider a $q$-deformed Clifford algebra $\A_q=\A_q(n)$ \cite{Ha}, which is the associative $K$-algebra with $1$ generated by $\psp_{a}$, $\psm_a$, $\om_a$, and $\om^{-1}_a$ for $a\in [\ov{n}]$
subject to the following relations:
{\allowdisplaybreaks
\begin{gather*}
\om_a \om_b= \om_b\om_a, \ \ \om_a\om_a^{-1}=1, \\
\om_a \psp_b \om_a^{-1} =
q^{ \delta_{ab}}\psp_b ,\ \ \ \ \ \om_a \psm_b \om_a^{-1} =
q^{- \delta_{ab}}\psm_b  \\
\psp_a\psp_b+ \psp_b\psp_a=0,\ \ \ \ \psm_a\psm_b+ \psm_b\psm_a=0 , \\
\psp_a\psm_b+ \psm_b\psp_a=0 \ \ \ \ (a\neq b), \\
\psp_a\psm_a =\frac{q\om_a-q^{-1}\om_a^{-1}}{q-q^{-1}},\ \ \
\psm_a\psp_a =-\frac{\om_a-\om_a^{-1}}{q-q^{-1}}.
\end{gather*}}

Let $\mathscr{E}_q$ be the left $\A_q$-module generated by $|0\rangle$ satisfying
$\psm_a|0\rangle =0$ and $\omega_a|0\rangle =q^{-1}|0\rangle$ 
for $a\in[\ov{n}]$. 
Then $\mathscr{E}_q$ is an irreducible $\A_q$-module with a $K$-linear basis $\{\,\psi_{{\bf m}}|0\rangle\,|\,{\bf m}\in \B\,\}$ (cf. \cite[Proposition 2.1]{Ha}), where ${\bf B}=\{\,(m_{a})\,|\,  a\in[\ov{n}],\  m_a\in\Z_2 \,\}$ and
$\psi_{{\bf m}}|0\rangle=\psp_{\ov{n}}^{m_{\ov{n}}}\cdots \psp_{\ov{1}}^{m_{\ov{1}}}|0\rangle$  for ${\bf m}=(m_a)\in \bf{B}$.

We put
\begin{equation}
\E=\mathscr{E}_{q_1},
\end{equation}
where $q_{_1}=q^{(\alpha_1,\alpha_1)/2}$  is equal to $q^{1/2}$ for  $C_n^{(1)}$, $q^2$ for $D_{n+1}^{(2)}$, and $q$ otherwise. One may regard $\E$ as an exterior algebra generated by an $n$-dimensional space $V$ with basis $\{v_{\ov{n}},\ldots, v_{\ov{1}}\}$ by identifying $\psp_{i_1}\cdots \psp_{i_k}|0\rangle$ with $v_{i_1}\wedge \cdots \wedge v_{i_k}$ for $\ov{n}\leq i_1<\ldots< i_k\leq \ov{1}$. Here, we understand $V$ as the dual natural representation of $U_q(\mf{sl}_n)\subset U'_q({\mf g})$. Then 

\begin{prop}[Theorem 3.2 in \cite{Ha}]\label{rho-1}
$\E$ has a $U_q(\mf{sl}_{n})$-module structure, where 
\begin{align*}
 t_i \longmapsto \omega_{\ov{i+1}}\omega_{\ov{i}}^{-1}, &\   \  \
 e_i \longmapsto \psp_{\ov{i+1}}\psm_{\ov{i}} ,  \ \ \  f_i \longmapsto \psm_{\ov{i+1}}\psp_{\ov{i}},
\end{align*}
for $i=1,\ldots,n-1$.
\end{prop}

Let $(\ ,\  )_{\E}$ be a non-degenerate symmetric bilinear form on $\E$ such that 
\begin{equation}\label{polarization on E}
(\psp_{\bf m}|0\rangle,\psp_{\bf m'}|0\rangle)_{\E}=\delta_{{\bf m}\,{\bf m}'},
\end{equation}
for ${\bf m},{\bf m}'\in {\bf B}$. Then it is straightforward to check that $\E$ has a polarizable crystal base $(L(\E),B(\E))$ with respect to $(\ , \ )_{\E}$, where
\begin{equation*}
L(\E)=\sum_{{\bf m}\in {\bf B}}\mathbb{A}\,\psp_{\bf m}|0\rangle, \ \ \  B(\E)=\{\,\psp_{\bf m}|0\rangle \!\!\!\pmod{qL(\E)} \,|\,{\bf m}\in {\bf B}\,\}.
\end{equation*}
We may identify $B(\E)$ with ${\bf B}$, and we have   for $i=1,\ldots,n-1$ and ${\bf m}=(m_a)\in {\bf B}$,  
\begin{equation}\label{crystal of type A}
\begin{split}
\te_i {\bf m} &=
\begin{cases}
{\bf m}+{\bf e}_{\ov{i+1}}-{\bf e}_{\ov{i}},& \text{if $(m_{\ov{i+1}},m_{\ov{i}})=(0,1)$},\\
0, & \text{otherwise},
\end{cases}\\
\tf_i {\bf m} &=
\begin{cases}
{\bf m}-{\bf e}_{\ov{i+1}}+{\bf e}_{\ov{i}},& \text{if $(m_{\ov{i+1}},m_{\ov{i}})=(1,0)$},\\
0, & \text{otherwise},
\end{cases}
\end{split}
\end{equation}
where ${\bf e}_a\in {\bf B}$ has $1$ at the $a$-th component and $0$ elsewhere for $a\in [\ov{n}]$.

\subsection{$U'_q(\mf{g})$-module structure on $\E^{\otimes 2}$} Now, we will construct a $U'_q(\mf{g})$-module structure on $\E$ or $\E^{\otimes 2}$ by extending the action of  $U_q(\mf{sl}_{n})$.
\begin{prop}\label{local-1}
Suppose that $\diamondsuit_r=
\ {\def\lr#1{\multicolumn{1}{|@{\hspace{.3ex}}c@{\hspace{.3ex}}|}{\raisebox{-.5ex}{$#1$}}}\raisebox{-.2ex}
{$\scalebox{0.5}{\begin{array}[b]{c}
\cline{1-1}
\lr{\ \ \ }\\ 
\cline{1-1}
\end{array}}$}}$\, or \, $ 
 {\def\lr#1{\multicolumn{1}{|@{\hspace{.4ex}}c@{\hspace{.4ex}}|}{\raisebox{-.65ex}{$#1$}}}\raisebox{-1ex}
{$\scalebox{0.5}{\begin{array}[b]{c}
\cline{1-1}
\lr{\ \ \ }\\ 
\cline{1-1}
\lr{\ \ \ }\\
\cline{1-1}
\end{array}}$}} $ \, for some $r\in\{0,n\}$. Then $\E$ has a $U_q(\mf{g}_{n-r})$-module structure, where the action of $U_q(\mf{g}_0\cap\mf{g}_n)$ is as in Proposition \ref{rho-1} and 
\begin{equation*}
\begin{split}
\begin{cases}
t_0  \longmapsto  q_{_0}  \omega_{\ov{1}}, \ \ \ \ \ e_0 \longmapsto \psp_{\ov{1}}, \ \ \ \ \ f_0 \longmapsto  \psm_{\ov{1}}, & 
\text{if $\diamondsuit_0=\, {\def\lr#1{\multicolumn{1}{|@{\hspace{.6ex}}c@{\hspace{.6ex}}|}{\raisebox{-.3ex}{$#1$}}}\raisebox{-.2ex}
{$\scalebox{0.45}{\begin{array}[b]{c}
\cline{1-1}
\lr{\ \ \ }\\ 
\cline{1-1}
\end{array}}$}}$}\, ,\\
t_0  \longmapsto  q_{_0} \omega_{\ov{1}}  \omega_{\ov{2}}, \ \ e_0 \longmapsto \psp_{\ov{1}}\psp_{\ov{2}}, \ \ f_0 \longmapsto \psm_{\ov{2}}\psm_{\ov{1}}, & 
\text{if $\diamondsuit_0=\,  {\def\lr#1{\multicolumn{1}{|@{\hspace{.6ex}}c@{\hspace{.6ex}}|}{\raisebox{-.3ex}{$#1$}}}\raisebox{-1.4ex}
{$\scalebox{0.45}{\begin{array}[b]{c}
\cline{1-1}
\lr{\ \ \ }\\ 
\cline{1-1}
\lr{\ \ \ }\\
\cline{1-1}
\end{array}}$}}$}\, ,\\
\end{cases}
\end{split}
\end{equation*}
\begin{equation*}
\begin{split}
\begin{cases}
t_n  \longmapsto   q_n^{-1} \omega_{\ov{n}}^{-1}, \ \ \ \ \ \ \ \ \ \ \ e_n \longmapsto \psm_{\ov{n}}, \ \ \ \ \ \ \ \ f_n \longmapsto  \psp_{\ov{n}}, & 
\text{if $\diamondsuit_n=\ {\def\lr#1{\multicolumn{1}{|@{\hspace{.6ex}}c@{\hspace{.6ex}}|}{\raisebox{-.3ex}{$#1$}}}\raisebox{-.3ex}
{$\scalebox{0.45}{\begin{array}[b]{c}
\cline{1-1}
\lr{\ \ \ }\\ 
\cline{1-1}
\end{array}}$}}$}\, ,\\
t_n  \longmapsto  q_n^{-1} (\omega_{\ov{n}}  \omega_{\ov{n-1}})^{-1}, \ \ e_n \longmapsto \psm_{\ov{n}}\psm_{\ov{n-1}}, \ \ f_n \longmapsto \psp_{\ov{n-1}}\psp_{\ov{n}}, & 
\text{if $\diamondsuit_n=\  {\def\lr#1{\multicolumn{1}{|@{\hspace{.6ex}}c@{\hspace{.6ex}}|}{\raisebox{-.3ex}{$#1$}}}\raisebox{-1.4ex}
{$\scalebox{0.45}{\begin{array}[b]{c}
\cline{1-1}
\lr{\ \ \ }\\ 
\cline{1-1}
\lr{\ \ \ }\\
\cline{1-1}
\end{array}}$}}$}\, ,
\end{cases}
\end{split}
\end{equation*}
and $(L(\E),B(\E))$ is a polarizable crystal base of $\E$ as  a $U_q(\mf{g}_{n-r})$-module with respect to $(\ , \ )_{\E}$.
\end{prop} 
\pf Suppose that $r=n$. Then $\E$ is a $U_q(\mf{g}_0)$-module by \cite[Theorem 4.1]{Ha} with a little modification (cf. \cite[Proposition 5.3]{K13-2} on which our presentation is based on), and $(L(\E),B(\E))$ is its crystal base as a $U_q(\mf{g}_0)$-module by \cite[Theorem 5.6]{K13-2}. It is also easy to check that $(L(\E),B(\E))$ is polarizable. The proof for $r=0$ is almost the same. 
\qed\vskip 2mm

Under the hypothesis of Proposition \ref{local-1}, we have for ${\bf m}=(m_a)\in {\bf B}$,
\begin{equation}\label{crystal of type B,D}
\begin{split}
\te_r {\bf m}&=
\begin{cases}
{\bf m}-{\bf e}_{\ov{n}},& \text{if $r=n$ with $\diamondsuit_n=\, {\def\lr#1{\multicolumn{1}{|@{\hspace{.6ex}}c@{\hspace{.6ex}}|}{\raisebox{-.3ex}{$#1$}}}\raisebox{-.2ex}
{$\scalebox{0.45}{\begin{array}[b]{c}
\cline{1-1}
\lr{\ \ \ }\\ 
\cline{1-1}
\end{array}}$}}$\, and $m_{\ov{n}}=1$},\\
{\bf m}+{\bf e}_{\ov{1}},& \text{if $r=0$ with\, $\diamondsuit_0=\, {\def\lr#1{\multicolumn{1}{|@{\hspace{.6ex}}c@{\hspace{.6ex}}|}{\raisebox{-.3ex}{$#1$}}}\raisebox{-.2ex}
{$\scalebox{0.45}{\begin{array}[b]{c}
\cline{1-1}
\lr{\ \ \ }\\ 
\cline{1-1}
\end{array}}$}}$\, and $m_{\ov{1}}=0$},\\
{\bf m}-{\bf e}_{\ov{n}}-{\bf e}_{\ov{n-1}},& \text{if $r=n$  with $\diamondsuit_n= \,
 {\def\lr#1{\multicolumn{1}{|@{\hspace{.6ex}}c@{\hspace{.6ex}}|}{\raisebox{-.3ex}{$#1$}}}\raisebox{-1.2ex}
{$\scalebox{0.45}{\begin{array}[b]{c}
\cline{1-1}
\lr{\ \ \ }\\ 
\cline{1-1}
\lr{\ \ \ }\\
\cline{1-1}
\end{array}}$}}$\, and $m_{\ov{n}}=m_{\ov{n-1}}=1$},\\
{\bf m}+{\bf e}_{\ov{2}}+{\bf e}_{\ov{1}},& \text{if $r=0$  with\, $\diamondsuit_0= \,
 {\def\lr#1{\multicolumn{1}{|@{\hspace{.6ex}}c@{\hspace{.6ex}}|}{\raisebox{-.3ex}{$#1$}}}\raisebox{-1.4ex}
{$\scalebox{0.45}{\begin{array}[b]{c}
\cline{1-1}
\lr{\ \ \ }\\ 
\cline{1-1}
\lr{\ \ \ }\\
\cline{1-1}
\end{array}}$}}$\, and $m_{\ov{2}}=m_{\ov{1}}=0$},\\
0, & \text{otherwise}.
\end{cases}
\end{split}
\end{equation}
Recall that $\tf_r{\bf m}={\bf m}'$ is determined by the relation $\te_r {\bf m}'={\bf m}$. In fact, $\B$ is the crystal of the spin representation (resp. the sum of two spin representations)
when  $\mf{g}_r=B_n$ (resp. $D_n$) \cite{KN}.

\begin{prop}\label{local-2}
Suppose that $\diamondsuit_r= {\def\lr#1{\multicolumn{1}{|@{\hspace{.3ex}}c@{\hspace{.3ex}}|}{\raisebox{-.3ex}{$#1$}}}\raisebox{-.2ex}
{$\scalebox{0.5}{\begin{array}[b]{cc}
\cline{1-1}\cline{2-2}
\lr{\ \ \ }&\lr{\ \ \ }\\ 
\cline{1-1}\cline{2-2}
\end{array}}$}}$\,  for some $r\in\{0,n\}$. Then $\E^{\otimes 2}$ has a $U_q(\mf{g}_{n-r})$-module structure, where the action of $U_q(\mf{g}_0\cap\mf{g}_n)$ is as in Proposition \ref{rho-1} and   
\begin{equation*}
\begin{split}
\begin{cases}
t_0  \longmapsto   q_{_0} \omega_{\ov{1}}\otimes  \omega_{\ov{1}}, \ \ \ \ \ \ \ e_0 \longmapsto \psp_{\ov{1}}\otimes  \psp_{\ov{1}}, \ \ \ \ \ f_0 \longmapsto  \psm_{\ov{1}}\otimes  \psm_{\ov{1}},  \\
t_n  \longmapsto  q_{n}^{-1} \omega_{\ov{n}}^{-1} \otimes  \omega_{\ov{n}}^{-1}, \ \, e_n \longmapsto \psm_{\ov{n}}\otimes  \psm_{\ov{n}}, \ \ \ \ f_n \longmapsto \psp_{\ov{n}}\otimes  \psp_{\ov{n}},
\end{cases}
\end{split}
\end{equation*}
and $(L(\E)^{\otimes 2},B(\E)^{\otimes 2})$ is a polarizable crystal base of $\E^{\otimes 2}$ as  a $U_q(\mf{g}_{n-r})$-module with respect to $(\ , \ )_{\E^{\otimes 2}}$, which is induced from $(\ , \ )_{\E}$.
\end{prop} 
\pf The proof is similar to that of Proposition \ref{local-1}.\qed \vskip 2mm

Under the hypothesis of Proposition \ref{local-2}, we have for ${\bf m}\otimes {\bf m}'=(m_a)\otimes (m'_a)\in {\bf B}^{\otimes 2}$
\begin{equation}\label{crystal of type C}
\begin{split}
\te_r ({\bf m}\otimes {\bf m}')&=
\begin{cases}
({\bf m}-{\bf e}_{\ov{n}})\otimes ({\bf m}'-{\bf e}_{\ov{n}}),& \text{if $r=n$ and $m_{\ov{n}}=m'_{\ov{n}}=1$},\\
({\bf m}+{\bf e}_{\ov{1}})\otimes ({\bf m}'+{\bf e}_{\ov{1}}),& \text{if $r=0$ and $m_{\ov{1}}=m'_{\ov{1}}=0$},\\
0, & \text{otherwise}.
\end{cases}
\end{split}
\end{equation}

Now, we have the following.

\begin{prop}\label{main-1} Let ${\mf g}$ be an affine Kac-Moody algebra of type $\pmb{\diamondsuit}=(\diamondsuit_0, \diamondsuit_n)$.
\begin{itemize}
\item[(1)] If $\diamondsuit_0, \diamondsuit_n \neq {\def\lr#1{\multicolumn{1}{|@{\hspace{.3ex}}c@{\hspace{.3ex}}|}{\raisebox{-.3ex}{$#1$}}}\raisebox{-.2ex}
{$\scalebox{0.5}{\begin{array}[b]{cc}
\cline{1-1}\cline{2-2}
\lr{\ \ \ }&\lr{\ \ \ }\\ 
\cline{1-1}\cline{2-2}
\end{array}}$}}$\ , then $\E$ is a finite-dimensional semisimple $U'_q(\mf{g})$-module  with a polarizable crystal base $(L(\E),B(\E))$ with ${\rm wt}(|0\rangle)={\rm cl}(\varpi_n)$. 

\item[(2)] $\E^{\otimes 2}$ is a  finite-dimensional semisimple $U'_q(\mf{g})$-module  with a polarizable crystal base $(L(\E)^{\otimes 2},B(\E)^{\otimes 2})$ with
\begin{equation*}
{\rm wt}(|0\rangle\otimes |0\rangle)=
\begin{cases}
2{\rm cl}(\varpi_n), & \text{if $\diamondsuit_0, \diamondsuit_n \neq {\def\lr#1{\multicolumn{1}{|@{\hspace{.6ex}}c@{\hspace{.6ex}}|}{\raisebox{-.08ex}{$#1$}}}\raisebox{-.1ex}
{$\scalebox{0.4}{\begin{array}[b]{cc}
\cline{1-1}\cline{2-2}
\lr{\ \ \ }&\lr{\ \ \ }\\ 
\cline{1-1}\cline{2-2}
\end{array}}$}}$\ },\\
{\rm cl}(\varpi_n), & \text{if $\diamondsuit_0$ or $\diamondsuit_n ={\def\lr#1{\multicolumn{1}{|@{\hspace{.6ex}}c@{\hspace{.6ex}}|}{\raisebox{-.08ex}{$#1$}}}\raisebox{-.1ex}
{$\scalebox{0.4}{\begin{array}[b]{cc}
\cline{1-1}\cline{2-2}
\lr{\ \ \ }&\lr{\ \ \ }\\ 
\cline{1-1}\cline{2-2}
\end{array}}$}}$\ }.
\end{cases}
\end{equation*}
\end{itemize}
\end{prop}
\pf  It follows from Propositions \ref{local-1} and \ref{local-2} that $\E^{\otimes N}$ ($N=1,2$) is a $U_q'(\mf{g})$-module, and $(L(\E)^{\otimes N},B(\E)^{\otimes N})$ is its polarizable crystal base, which also implies that $\E^{\otimes N}$ is semisimple by Proposition \ref{polarizable implies semisimple}. 
\qed 

\subsection{Binary matrices and crystal of $\E^{\otimes 2}$}\label{crystal structure on M}
Let ${\bf M}$ be the set of binary matrices ${\bf m}=[\,m_{ab}\,]$ ($a\in[\ov{n}],b\in \{1,2\}$). For ${\bf m}\in {\bf M}$, let ${\bf m}_{(a)}=[\,m_{a1}\  m_{a2}\,]$ be the $a$-th row  and ${\bf m}^{(b)}=[\,m_{ab}\,]$ the $b$-th column of ${\bf m}$.  We put $|{\bf m}^{(b)}|=\sum_{a}m_{ab}$ for $b=1,2$, and $|{\bf m}|=|{\bf m}^{(1)}|+|{\bf m}^{(2)}|$.

By Proposition \ref{main-1} (2), we may regard ${\bf M}$ as a crystal of $\E^{\otimes 2}$ identifying ${\bf m}\in {\bf M}$ with $\psp_{{\bf m}^{(1)}}|0\rangle \otimes \psp_{{\bf m}^{(2)}}|0\rangle \in B(\E)^{\otimes 2}={\bf B}^{\otimes 2}$. 

\begin{figure}
\includegraphics[width=9.5cm, height=20cm]{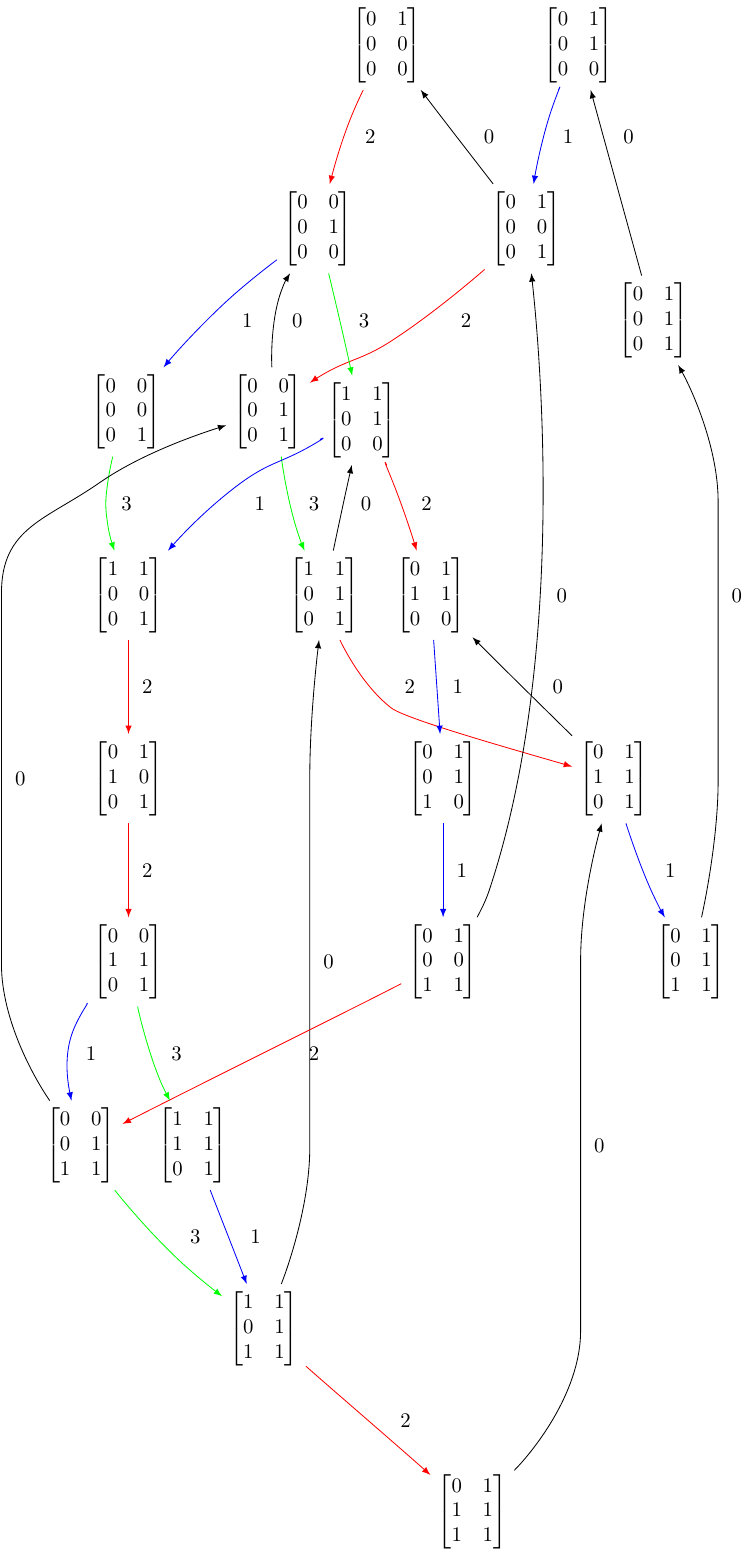}
\caption{ $B(W(\varpi_2))=C({\bf v}_{2,0})$ of type $A_6^{(2)}$.} \label{Graph B}
\end{figure}

Let us describe $\te_r$ for $r\in \{0,n\}$ on ${\bf M}$ using \eqref{crystal of type B,D}, \eqref{crystal of type C}, and the tensor product rule of crystals (cf. \cite{HK}). This will be useful for the arguments in the next sections (see  Figures \ref{Graph A} and \ref{Graph B},  for example).
\vskip 2mm

Let ${\bf m}\in {\bf M}$ be given.\vskip 2mm

 \textsc{Case 1}. $\diamondsuit_r=\ {\def\lr#1{\multicolumn{1}{|@{\hspace{.4ex}}c@{\hspace{.4ex}}|}{\raisebox{-.6ex}{$#1$}}}\raisebox{-.3ex}
{$\scalebox{0.5}{\begin{array}[b]{c}
\cline{1-1}
\lr{\ \ \ }\\ 
\cline{1-1}
\end{array}}$}}$\ .

Since $\te_n {\bf m}$ and $\te_0 {\bf m}$  depend only on ${\bf m}_{(\ov{n})}$ and ${\bf m}_{(\ov{1})}$, respectively, it is enough to describe them in terms of ${\bf m}_{(\ov{n})}$ and ${\bf m}_{(\ov{1})}$. We have  
\begin{equation}\label{Kashiwara op for (1)}
\te_n {\bf m}_{(\ov{n})}=
\begin{cases}
[\, 0\ 0\, ], & \text{if ${\bf m}_{(\ov{n})}=[\,1\ 0\,]$},  \\
[\,1\ 0\,], & \text{if ${\bf m}_{(\ov{n})}=[\,1\ 1\,]$}, \\
0, &  \text{otherwise},
\end{cases}\ \ \ 
\te_0 {\bf m}_{(\ov{1})}=
\begin{cases}
[\,0\ 1\,], & \text{if ${\bf m}_{(\ov{1})}=[\,0\ 0\,]$},  \\
[\,1\ 1\,], & \text{if ${\bf m}_{(\ov{1})}=[\,0\ 1\,]$}, \\
0, &  \text{otherwise}.
\end{cases}
\end{equation}

 \textsc{Case 2}. Suppose that $\diamondsuit_r={\def\lr#1{\multicolumn{1}{|@{\hspace{.3ex}}c@{\hspace{.3ex}}|}{\raisebox{-.3ex}{$#1$}}}\raisebox{-.2ex}
{$\scalebox{0.5}{\begin{array}[b]{cc}
\cline{1-1}\cline{2-2}
\lr{\ \ \ }&\lr{\ \ \ }\\ 
\cline{1-1}\cline{2-2}
\end{array}}$}}$\ .

As in Case 1, $\te_n {\bf m}$ and $\te_0 {\bf m}$  depend only on ${\bf m}_{(\ov{n})}$ and ${\bf m}_{(\ov{1})}$, respectively. We have
\begin{equation}\label{Kashiwara op for (2)}
\te_n {\bf m}_{(\ov{n})}=
\begin{cases}
[\,0\ 0\,], & \text{if ${\bf m}_{(\ov{n})}=[\,1\ 1\,]$},  \\
0, &  \text{otherwise},
\end{cases}
%\ \ \ \tf_n {\bf m}_{(\ov{n})}=
%\begin{cases}
%[1\ 1], & \text{if ${\bf m}_{(\ov{n})}=[0\ 0]$},  \\
%0, &  \text{otherwise},
%\end{cases}
\ \ \
\te_0 {\bf m}_{(\ov{1})}=
\begin{cases}
[\,1\ 1\,], & \text{if ${\bf m}_{(\ov{1})}=[\,0\ 0\,]$},  \\
0, &  \text{otherwise}.
\end{cases}
%\ \ \ \tf_0 {\bf m}_{(\ov{1})}=
%\begin{cases}
%[0\ 0], & \text{if ${\bf m}_{(\ov{1})}=[1\ 1]$},  \\
%0, &  \text{otherwise}.
%\end{cases}
\end{equation}

 \textsc{Case 3}. Suppose that $\diamondsuit_r=\ {\def\lr#1{\multicolumn{1}{|@{\hspace{.6ex}}c@{\hspace{.6ex}}|}{\raisebox{-.3ex}{$#1$}}}\raisebox{-.9ex}
{$\scalebox{0.5}{\begin{array}[b]{c}
\cline{1-1}
\lr{\ \ \ }\\ 
\cline{1-1}
\lr{\ \ \ }\\
\cline{1-1}
\end{array}}$}} $\ .\vskip 2mm

It is enough to describe in terms of  the $2\times 2$-submatrix
\renewcommand\arraystretch{0.9}$
\begin{bmatrix}
{\bf m}_{(\ov{2})} \\ {\bf m}_{(\ov{1})} 
\end{bmatrix}$ or
$\renewcommand\arraystretch{0.97}
\begin{bmatrix}
{\bf m}_{(\ov{n})} \\ {\bf m}_{(\ov{n-1})} 
\end{bmatrix}=\renewcommand\arraystretch{0.9}
\begin{bmatrix}
\,p &  q\, \\
\,s  & t\,
\end{bmatrix}
$. We have for $r=n$  
\begin{equation}\label{Kashiwara op for (1,1)-1}
\te_n  
\begin{bmatrix}
{\bf m}_{(\ov{n})} \\ {\bf m}_{(\ov{n-1})} 
\end{bmatrix}
=
\begin{cases}
\renewcommand\arraystretch{0.9}
\begin{bmatrix}
\,p &  0\, \\
\,s  &  0\,
\end{bmatrix}, 
& \text{if $
\renewcommand\arraystretch{0.9}\begin{bmatrix}
\,p \, \\ s 
\end{bmatrix}
\neq 
\begin{bmatrix}
\,0\,  \\ 0 
\end{bmatrix}$ and  $
\renewcommand\arraystretch{0.9}\begin{bmatrix}
\,q\,  \\ t
\end{bmatrix}
=
\begin{bmatrix}
\,1\,  \\ 1 
\end{bmatrix}
$},  \\  \\
\renewcommand\arraystretch{0.9}\begin{bmatrix}
\,0 &  q\, \\
\,0  &  t\,
\end{bmatrix}, 
& \text{if $
\renewcommand\arraystretch{0.9}\begin{bmatrix}
\,p\,  \\ s 
\end{bmatrix}
= 
\begin{bmatrix}
\,1\, \\ 1 
\end{bmatrix}$ and  $
\renewcommand\arraystretch{0.9}\begin{bmatrix}
\,q\,  \\
t
\end{bmatrix}
\neq
\begin{bmatrix}
\,1\,   \\ 1 
\end{bmatrix}
$}, \\ \\ 
\ \ 0, &  \text{otherwise},\\
\end{cases}
\end{equation}
and for $r=0$
\begin{equation}\label{Kashiwara op for (1,1)-2}
\te_0  
\begin{bmatrix}
{\bf m}_{(\ov{2})} \\ {\bf m}_{(\ov{1})} 
\end{bmatrix}
=
\begin{cases}
\renewcommand\arraystretch{0.9}
\begin{bmatrix}
\,p & 1\, \\
\,s &  1\,
\end{bmatrix}, 
& \text{if  $\renewcommand\arraystretch{0.9}
\begin{bmatrix}
\,p\,  \\
s
\end{bmatrix}
\neq \renewcommand\arraystretch{0.9}
\begin{bmatrix}
\,1\,   \\ 1 
\end{bmatrix}
$ and 
$\renewcommand\arraystretch{0.9}
\begin{bmatrix}
\,q\,  \\
t
\end{bmatrix}
=\renewcommand\arraystretch{0.9}
\begin{bmatrix}
\,0\,   \\ 0 
\end{bmatrix}$},  \\  \\
\renewcommand\arraystretch{0.9}\begin{bmatrix}
\,1 & q \,\\
\,1 & t \,
\end{bmatrix}, 
& \text{if 
$\renewcommand\arraystretch{0.9}
\begin{bmatrix}
\,p \,  \\
s
\end{bmatrix}
=
\begin{bmatrix}
\,0 \,   \\ 0 
\end{bmatrix}
$ and
${\renewcommand\arraystretch{0.9}
\begin{bmatrix}
\,q\,  \\
t
\end{bmatrix}
\neq 
\begin{bmatrix}
\,0\, \\ 0 
\end{bmatrix}}$}, \\ \\ 
\ \ 0, &  \text{otherwise}.
\end{cases}
\end{equation}\vskip 2mm

There is also a $U_q(\mf{sl}_2)$-crystal structure on ${\bf M}$.
 For each $a\in [\ov{n}]$, we may regard ${\bf m}_{(a)}$ as a crystal element over $U_q(\mf{sl}_2)$ with Kashiwara operators $\td{{\it E}}$ and $\td{{\it F}}$ such that $\td{E}[\,0\ 1\,]=[\,1 \ 0\,]$, $\td{F}[\,1\ 0\,]=[\,0 \ 1\,]$, and $\td{X}[\,0\ 0\,]=\td{X}[\,1\ 1\,]=0$ ($X=E, F$), and understand ${\bf m}$ as ${\bf m}_{(\ov{1})}\otimes \cdots\otimes {\bf m}_{(\ov{n})}$. Then ${\bf M}$ is a 
$(U_q(\mf{sl}_{n}),U_q(\mf{sl}_2))$-bicrystal  (cf. \cite{DK,La}).
For ${\bf m}\in {\bf M}$, we put
\begin{equation}\label{sigma of m}
\sigma({\bf m})=(\varepsilon({\bf m}),\varphi({\bf m})),
\end{equation}
where $\varepsilon({\bf m})=\max\{\,k\,|\,\td{E}^k{\bf m}\neq 0\,\}$ and $\varphi({\bf m})=\max\{\,k\,|\,\td{F}^k{\bf m}\neq 0\,\}$.  Note that  if $\tilde{x}_i {\bf m}\neq 0$ for some $i\in I\setminus\{0,n\}$ and $x\in \{e,f\}$, then 
\begin{equation}\label{invariance of sigma}
\sigma(\tilde{x}_i {\bf m})=\sigma({\bf m}).
\end{equation}

\section{Crystal structure on $\E^{\otimes 2}$}

\subsection{Decomposition of the crystal of $\E^{\otimes 2}$} 
Suppose that ${\mf g}$ is of type $\pmb{\diamondsuit}=(\diamondsuit_0, \diamondsuit_n)$ with 
$\diamondsuit_0, \diamondsuit_n \neq{\def\lr#1{\multicolumn{1}{|@{\hspace{.3ex}}c@{\hspace{.3ex}}|}{\raisebox{-.3ex}{$#1$}}}\raisebox{-.2ex}
{$\scalebox{0.5}{\begin{array}[b]{cc}
\cline{1-1}\cline{2-2}
\lr{\ \ \ }&\lr{\ \ \ }\\ 
\cline{1-1}\cline{2-2}
\end{array}}$}}$\ .
Let ${\bf v}_n=(0,\ldots,0)$ and ${\bf v}_{n-1}={\bf e}_{\ov{n}}$. Then it is not difficult to see that
\begin{equation}\label{decomposition of B}
{\bf B}=
\begin{cases}
C({\bf v}_n)\sqcup C({\bf v}_{n-1}), & \text{if \ $\pmb{\diamondsuit}= (\, {\def\lr#1{\multicolumn{1}{|@{\hspace{.6ex}}c@{\hspace{.6ex}}|}{\raisebox{-.3ex}{$#1$}}}\raisebox{-.8ex}
{$\scalebox{0.35}{\begin{array}[b]{c}
\cline{1-1}
\lr{\ \ \ }\\ 
\cline{1-1}
\lr{\ \ \ }\\
\cline{1-1}
\end{array}}$}}\, ,\, {\def\lr#1{\multicolumn{1}{|@{\hspace{.6ex}}c@{\hspace{.6ex}}|}{\raisebox{-.3ex}{$#1$}}}\raisebox{-.8ex}
{$\scalebox{0.35}{\begin{array}[b]{c}
\cline{1-1}
\lr{\ \ \ }\\ 
\cline{1-1}
\lr{\ \ \ }\\
\cline{1-1}
\end{array}}$}} \, )$},\\
C({\bf v}_n), & \text{otherwise}, \\
\end{cases}
\end{equation}
where $C({\bf v})$ denotes the connected component of ${\bf v}$ in ${\bf B}$ as a $U'_q(\mf{g})$-crystal. 

Next, suppose that ${\mf g}$ is of type $\pmb{\diamondsuit}=(\diamondsuit_0, \diamondsuit_n)$ with 
  $\diamondsuit_0$ or $\diamondsuit_n ={\def\lr#1{\multicolumn{1}{|@{\hspace{.3ex}}c@{\hspace{.3ex}}|}{\raisebox{-.3ex}{$#1$}}}\raisebox{-.2ex}
{$\scalebox{0.5}{\begin{array}[b]{cc}
\cline{1-1}\cline{2-2}
\lr{\ \ \ }&\lr{\ \ \ }\\ 
\cline{1-1}\cline{2-2}
\end{array}}$}}$\ .
For $0\leq k\leq n$ and $0\leq l\leq n-k$, let ${\bf v}_{k,l}=\left[\, {\bf v}^{(1)}_{k,l}\  {\bf v}^{(2)}_{k,l}\,\right]\in {\bf M}$ be given by
\begin{equation*}
\begin{split}
{\bf v}^{(1)}_{k,l} = {\bf e}_{\ov{n}}+\cdots+{\bf e}_{\ov{n-l+1}},\ \ \
{\bf v}^{(2)}_{k,l} = {\bf e}_{\ov{n-l}}+\cdots+{\bf e}_{\ov{k+1}}.
\end{split}
\end{equation*}
Here we understand ${\bf e}_a$ as a column vector and ${\bf v}^{(1)}_{k,l}$ (resp. ${\bf v}^{(2)}_{k,l}$) as a zero vector when $l=0$ (resp. $l=n-k$). Note that the number of $1$'s in ${\bf v}_{k,l}$ is $n-k$, while the number of $1$'s in the first column ${\bf v}^{(1)}_{k,l}$ is $l$. 
We have $\te_i{\bf v}_{k,l}=0$ for all $i\in I\setminus\{0,n\}$ by \eqref{crystal of type A} and the tensor product rule of crystals, where 
\begin{equation}\label{weight of v_{k,l}}
{\rm wt}({\bf v}_{k,l})=
\begin{cases}
{\rm cl}(\varpi_k), & \text{if $1\leq k\leq n$},\\
0, & \text{if $k=0$}.
\end{cases}
\end{equation}

For ${\bf m}\in {\bf M}$, let $C({\bf m})$ denote the connected component of ${\bf m}$ in $\bf{M}$  as a $U'_q(\mf{g})$-crystal. Then we have the following decomposition of ${\bf M}$. 
\begin{prop}\label{decomposition of M}
 Suppose that ${\mf g}$ is of type $\pmb{\diamondsuit}=(\diamondsuit_0, \diamondsuit_n)$ with 
  $\diamondsuit_0$ or $\diamondsuit_n ={\def\lr#1{\multicolumn{1}{|@{\hspace{.3ex}}c@{\hspace{.3ex}}|}{\raisebox{-.3ex}{$#1$}}}\raisebox{-.2ex}
{$\scalebox{0.5}{\begin{array}[b]{cc}
\cline{1-1}\cline{2-2}
\lr{\ \ \ }&\lr{\ \ \ }\\ 
\cline{1-1}\cline{2-2}
\end{array}}$}}$\ .
Then as a $U'_q(\mf{g})$-crystal,  
\begin{equation*}
{\bf M}=\bigsqcup_{(k,l)\in  H^{\pmb{\diamondsuit}}}C({\bf v}_{k,l}),
\end{equation*}
where
\begin{equation*}
{H}^{\pmb{\diamondsuit}}=
\begin{cases}
\{\, (k,l)\,|\, 0\leq k\leq  n,\, 0\leq l \leq  n-k \,\}, & \text{if \ $\pmb{\diamondsuit}=(\, {\def\lr#1{\multicolumn{1}{|@{\hspace{.6ex}}c@{\hspace{.6ex}}|}{\raisebox{-.08ex}{$#1$}}}\raisebox{-.1ex}
{$\scalebox{0.35}{\begin{array}[b]{cc}
\cline{1-1}\cline{2-2}
\lr{\ \ \ }&\lr{\ \ \ }\\ 
\cline{1-1}\cline{2-2}
\end{array}}$}}\, ,\, 
{\def\lr#1{\multicolumn{1}{|@{\hspace{.6ex}}c@{\hspace{.6ex}}|}{\raisebox{-.08ex}{$#1$}}}\raisebox{-.1ex}
{$\scalebox{0.35}{\begin{array}[b]{cc}
\cline{1-1}\cline{2-2}
\lr{\ \ \ }&\lr{\ \ \ }\\ 
\cline{1-1}\cline{2-2}
\end{array}}$}}\,)$},\\
\{\, (k,n-k)\,|\, 0\leq k \leq n  \,\}, & \text{if \ $\pmb{\diamondsuit}=(\, {\def\lr#1{\multicolumn{1}{|@{\hspace{.6ex}}c@{\hspace{.6ex}}|}{\raisebox{-.3ex}{$#1$}}}\raisebox{-.1ex}
{$\scalebox{0.35}{\begin{array}[b]{c}
\cline{1-1}
\lr{\ \ \ }\\ 
\cline{1-1}
\end{array}}$}}\, ,\, {\def\lr#1{\multicolumn{1}{|@{\hspace{.6ex}}c@{\hspace{.6ex}}|}{\raisebox{-.08ex}{$#1$}}}\raisebox{-.1ex}
{$\scalebox{0.35}{\begin{array}[b]{cc}
\cline{1-1}\cline{2-2}
\lr{\ \ \ }&\lr{\ \ \ }\\ 
\cline{1-1}\cline{2-2}
\end{array}}$}} \, )$},\\
\{\, (k,0)\,|\, 0\leq k \leq n  \,\}, & \text{if \ $\pmb{\diamondsuit}=(\, {\def\lr#1{\multicolumn{1}{|@{\hspace{.6ex}}c@{\hspace{.6ex}}|}{\raisebox{-.08ex}{$#1$}}}\raisebox{-.1ex}
{$\scalebox{0.35}{\begin{array}[b]{cc}
\cline{1-1}\cline{2-2}
\lr{\ \ \ }&\lr{\ \ \ }\\ 
\cline{1-1}\cline{2-2}
\end{array}}$}} \,,\,{\def\lr#1{\multicolumn{1}{|@{\hspace{.6ex}}c@{\hspace{.6ex}}|}{\raisebox{-.3ex}{$#1$}}}\raisebox{-.1ex}
{$\scalebox{0.35}{\begin{array}[b]{c}
\cline{1-1}
\lr{\ \ \ }\\ 
\cline{1-1}
\end{array}}$}}\, )$},\\
\{\, (k,n-k)\,|\, 0\leq k \leq n  \,\}\cup\{(0,n-1)\}, & \text{if \ $\pmb{\diamondsuit}= (\, {\def\lr#1{\multicolumn{1}{|@{\hspace{.6ex}}c@{\hspace{.6ex}}|}{\raisebox{-.3ex}{$#1$}}}\raisebox{-.8ex}
{$\scalebox{0.35}{\begin{array}[b]{c}
\cline{1-1}
\lr{\ \ \ }\\ 
\cline{1-1}
\lr{\ \ \ }\\
\cline{1-1}
\end{array}}$}}\, ,\, {\def\lr#1{\multicolumn{1}{|@{\hspace{.6ex}}c@{\hspace{.6ex}}|}{\raisebox{-.08ex}{$#1$}}}\raisebox{-.1ex}
{$\scalebox{0.35}{\begin{array}[b]{cc}
\cline{1-1}\cline{2-2}
\lr{\ \ \ }&\lr{\ \ \ }\\ 
\cline{1-1}\cline{2-2}
\end{array}}$}} \, )$}.\\
\end{cases}
\end{equation*}
\end{prop}
\pf 
Let $C$ be a connected component in ${\bf M}$ with respect to $\te_i$ and $\tf_i$ for $i\in I$. Choose
 ${\bf m}\in C$ such that $|{\bf m}|$ is minimal. We may assume that ${\bf m}$ is of $\mf{g}_0$-highest weight since $|\te_i{\bf m}|\leq |{\bf m}|$ if $\te_i{\bf m}\neq 0$ for $i\in I_0$. 

\textsc{Case 1}. Suppose that $\pmb{\diamondsuit}=(\, {\def\lr#1{\multicolumn{1}{|@{\hspace{.3ex}}c@{\hspace{.3ex}}|}{\raisebox{-.3ex}{$#1$}}}\raisebox{-.2ex}
{$\scalebox{0.5}{\begin{array}[b]{cc}
\cline{1-1}\cline{2-2}
\lr{\ \ \ }&\lr{\ \ \ }\\ 
\cline{1-1}\cline{2-2}
\end{array}}$}}\, ,\, 
{\def\lr#1{\multicolumn{1}{|@{\hspace{.3ex}}c@{\hspace{.3ex}}|}{\raisebox{-.3ex}{$#1$}}}\raisebox{-.2ex}
{$\scalebox{0.5}{\begin{array}[b]{cc}
\cline{1-1}\cline{2-2}
\lr{\ \ \ }&\lr{\ \ \ }\\ 
\cline{1-1}\cline{2-2}
\end{array}}$}}\,)$.
We first note that ${\bf m}_{(\ov{n})}\neq [\,1\  1\,]$. Otherwise, $\te_{\ov{n}}{\bf m}\neq 0$.   Suppose that ${\bf m}_{(\ov{n})}=[\,0\  1\,]$. Then there exists $\ov{k+1}\in [\ov{n}]$ such that ${\bf m}_{\ov{k'}}=[\,0\  1\,]$ for $k+1\leq k'\leq n$ and ${\bf m}_{\ov{k'}}=[\,0\  0\,]$ otherwise, since $\te_i {\bf m}=0$ for $i\in I_0$. This implies that ${\bf m}={\bf v}_{k,0}$. Suppose that ${\bf m}_{(\ov{n})}=[\,1\  0\,]$. Let $\ov{k'+1}$ be the smallest such that ${\bf m}_{\ov{k'+1}}=[\,1\ 0\,]$. If ${\bf m}_{\ov{k'}}=[\,0 \ 0\,]$, then ${\bf m}_{\ov{k'}}=\cdots={\bf m}_{\ov{1}}=[\,0\ 0\,]$. If ${\bf m}_{\ov{k'}}=[\,0\ 1\,]$, then  as in the previous case, we have ${\bf m}_{\ov{k'}}=\cdots={\bf m}_{\ov{k+1}}=[\,0\ 1\,]$ and ${\bf m}_{\ov{k}}=\cdots={\bf m}_{\ov{1}}=[\,0\ 0\,]$ for some $k$. This implies that ${\bf m}={\bf v}_{k,l}$, where $l=n-k'$. Hence ${\bf m}={\bf v}_{k,l}$  and $C=C({\bf v}_{k,l})$ for some $0\leq k\leq n$ and $0\leq l\leq n-k$.

Let ${\bf m}'\in C({\bf v}_{k,l})$ be given. 
We see from \eqref{Kashiwara op for (2)} and \eqref{invariance of sigma} that $\sigma({\bf m}')=\sigma(\tilde{x}_i{\bf m}')$ for $i\in I$ and $x\in\{e,f\}$ with $\tilde{x}_i{\bf m}'\neq 0$, and hence $\sigma({\bf m}')=(n-k-l,l)=\sigma({\bf v}_{k,l})$.  Also by the signature rule of tensor product of crystals with respect to $U_q(\mf{sl}_2)$ (cf. \cite[Remark 2.1.2]{KN}), we have $|{\bf m}'|\geq n-k$.
This implies that $C({\bf v}_{k,l})=C({\bf v}_{k'',l''})$ if and only if $k=k''$ and $l=l''$ for $0\leq k,k''\leq n$, $0\leq l\leq n-k$, and $0\leq l''\leq n-k''$.  Hence we obtain the decomposition of $\bf{M}$. 

\textsc{Case 2}. Suppose that
 $\pmb{\diamondsuit}=(\, {\def\lr#1{\multicolumn{1}{|@{\hspace{.3ex}}c@{\hspace{.3ex}}|}{\raisebox{-.3ex}{$#1$}}}\raisebox{-.2ex}
{$\scalebox{0.5}{\begin{array}[b]{c}
\cline{1-1}
\lr{\ \ \ }\\ 
\cline{1-1}
\end{array}}$}}\, ,\,{\def\lr#1{\multicolumn{1}{|@{\hspace{.3ex}}c@{\hspace{.3ex}}|}{\raisebox{-.3ex}{$#1$}}}\raisebox{-.2ex}
{$\scalebox{0.5}{\begin{array}[b]{cc}
\cline{1-1}\cline{2-2}
\lr{\ \ \ }&\lr{\ \ \ }\\ 
\cline{1-1}\cline{2-2}
\end{array}}$}} \, )$. As in Case 1, we have ${\bf m}={\bf v}_{k,n-k'}$ for some $k, k'$ with $0\leq k\leq k'\leq n$. But if $k< k'$, then ${\bf m}$ is connected to ${\bf v}_{k',n-k'}$ by applying $\tf_i$'s for $i\in \{\, k'-1,\ldots, 1, 0\,\}$ (see \eqref{Kashiwara op for (1)}), which contradicts the minimality of $|{\bf m}|$. Hence, ${\bf m}={\bf v}_{k,n-k}$, and $C=C({\bf v}_{k,n-k})$. 

Let ${\bf m}'\in C({\bf v}_{k,n-k})$ be given. Suppose that  $|{\bf m'}^{(1)}|< n-k$. We may assume that ${\bf m}'$ is of $\g_0$-highest weight by  \eqref{Kashiwara op for (2)}.  Then by the same argument as in the previous paragraph, we see that ${\bf m}'$ is connected to ${\bf v}_{k',n-k'}$ with $n-k'<n-k$, which contradicts the minimality of $|{\bf m}|$.  Hence  $|{\bf m'}^{(1)}|\geq  n-k$. This implies that $C({\bf v}_{k,n-k})=C({\bf v}_{k',n-k'})$ if and only if $k=k'$ for $0\leq k, k'\leq n$. The proof for $\pmb{\diamondsuit}=(\, {\def\lr#1{\multicolumn{1}{|@{\hspace{.4ex}}c@{\hspace{.4ex}}|}{\raisebox{-.3ex}{$#1$}}}\raisebox{-.2ex}
{$\scalebox{0.5}{\begin{array}[b]{cc}
\cline{1-1}\cline{2-2}
\lr{\ \ \ }&\lr{\ \ \ }\\ 
\cline{1-1}\cline{2-2}
\end{array}}$}} \,,\,{\def\lr#1{\multicolumn{1}{|@{\hspace{.3ex}}c@{\hspace{.3ex}}|}{\raisebox{-.3ex}{$#1$}}}\raisebox{-.2ex}
{$\scalebox{0.5}{\begin{array}[b]{c}
\cline{1-1}
\lr{\ \ \ }\\ 
\cline{1-1}
\end{array}}$}}\, )$ is almost the same.

\textsc{Case 3.}  Suppose that $\pmb{\diamondsuit}=(\, {\def\lr#1{\multicolumn{1}{|@{\hspace{.6ex}}c@{\hspace{.6ex}}|}{\raisebox{-.3ex}{$#1$}}}\raisebox{-.7ex}
{$\scalebox{0.45}{\begin{array}[b]{c}
\cline{1-1}
\lr{\ \ \ }\\ 
\cline{1-1}
\lr{\ \ \ }\\
\cline{1-1}
\end{array}}$}}\, ,\, {\def\lr#1{\multicolumn{1}{|@{\hspace{.4ex}}c@{\hspace{.4ex}}|}{\raisebox{-.3ex}{$#1$}}}\raisebox{-.2ex}
{$\scalebox{0.5}{\begin{array}[b]{cc}
\cline{1-1}\cline{2-2}
\lr{\ \ \ }&\lr{\ \ \ }\\ 
\cline{1-1}\cline{2-2}
\end{array}}$}} \,  )$. Then  we have ${\bf m}={\bf v}_{k,n-k'}$ for some $k, k'$ with $0\leq k\leq k'\leq n$. If $k=n$, then ${\bf m}={\bf v}_{n,0}$. If $k=0$, then ${\bf m}={\bf v}_{0,n}$ or ${\bf v}_{0,n-1}$ since $\tf_{0}{\bf v}_{0,n-k'}\neq 0$ for $k'>1$ by \eqref{Kashiwara op for (1,1)-2}, which  contradicts the minimality of $|{\bf m}|$. Note that $C({\bf v}_{0,n})=\{{\bf v}_{0,n}\}$ and $C({\bf v}_{0,n-1})=\{{\bf v}_{0,n-1}\}$.

\begin{figure}
\includegraphics[width=8cm, height=20cm]{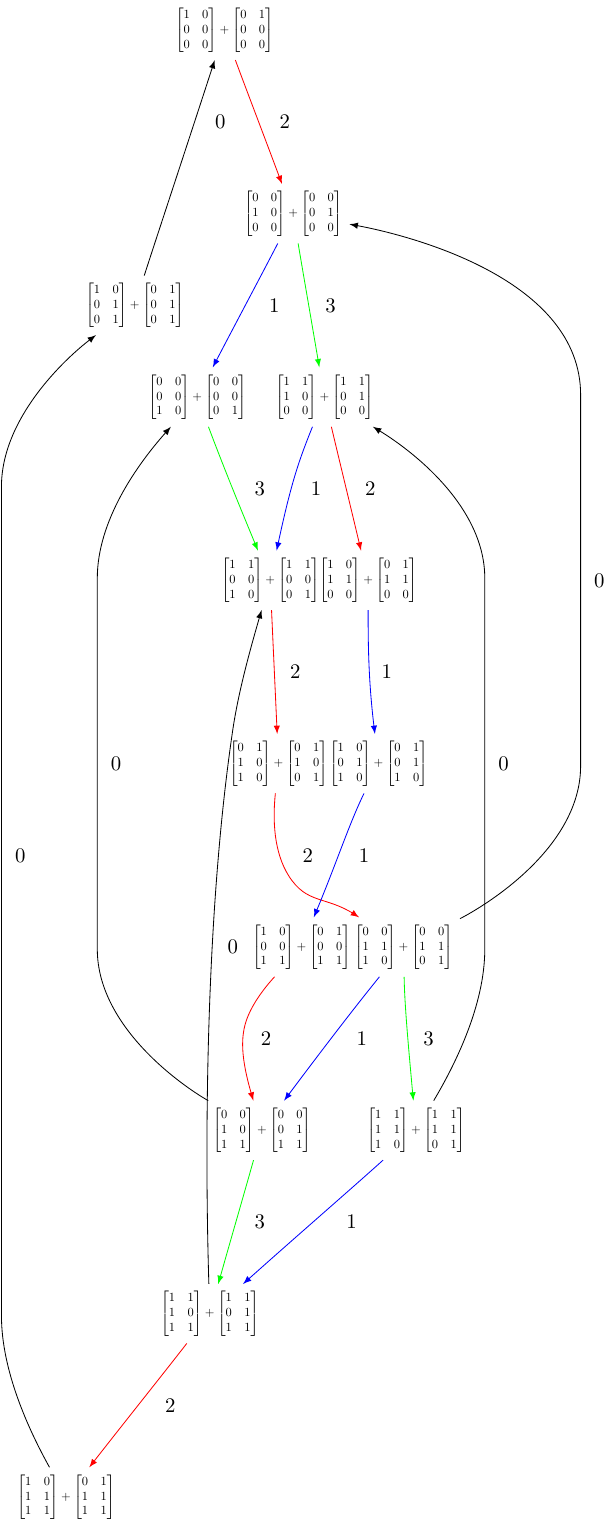}
\caption{ $B(W(\varpi_2))$ of type $A_5^{(2)}$ with vertices ${\bf m}+\varsigma({\bf m})$ for ${\bf m}\in C({\bf v}_{2,1})^+$. } \label{Graph C}
\end{figure}

If $1\leq k\leq n-1$, then ${\bf m}$ is connected to ${\bf v}_{k',n-k'}$ or ${\bf v}_{k'-1,n-k'}$ by applying $\tf_i$'s for $i\in \{\,k'-1,\ldots, 1, 0\,\}$. So, by the minimality of $|{\bf m}|$, we must have ${\bf m}={\bf v}_{k,n-k}$ or ${\bf v}_{k,n-k-1}$. On the other hand, let $w\in W$ be such that $w(\varpi_{k})=\varpi_{k}+\delta$. For example, put
\begin{equation}\label{delta shifting w}
w=(w_{k}w_{k-1}\cdots w_{1})s_n s_{n-1} \cdots s_{2}s_0 (w_{k}w_{k-1}\cdots w_{2})^{-1},
\end{equation}
where $w_i=s_is_{i+1}\cdots s_{n-k+i-1}$ for $i\in \{1,\ldots,k\}$. 
Then it is straightforward to check that 
\begin{equation}\label{delta shifting w-2}
{\texttt S}_w {\bf v}_{k,n-k-1}={\bf v}_{k,n-k}, \ \ {\texttt S}_w {\bf v}_{k,n-k}={\bf v}_{k,n-k-1},
\end{equation} 
(see Figure \ref{Graph C}).  This implies that $C=C({\bf v}_{k,n-k})=C({\bf v}_{k,n-k-1})$. 

Finally, by similar arguments as in Case 2, we can check that  $C({\bf v}_{k,n-k})=C({\bf v}_{k',n-k'})$ if and only if $k=k'$ for $1\leq k, k'\leq n$. Hence we have the decomposition of ${\bf M}$.
\qed

\subsection{Decomposition into classical crystals}  Suppose that ${\mf g}$ is of type $\pmb{\diamondsuit}=(\diamondsuit_0, \diamondsuit_n)$ with 
  $\diamondsuit_0$ or $\diamondsuit_n ={\def\lr#1{\multicolumn{1}{|@{\hspace{.4ex}}c@{\hspace{.4ex}}|}{\raisebox{-.3ex}{$#1$}}}\raisebox{-.2ex}
{$\scalebox{0.5}{\begin{array}[b]{cc}
\cline{1-1}\cline{2-2}
\lr{\ \ \ }&\lr{\ \ \ }\\ 
\cline{1-1}\cline{2-2}
\end{array}}$}}$\ . For a $\mf{g}_0$-dominant weight  $\lambda\in P$, let $B_0(\lambda)$ be the crystal of the irreducible $U_q(\mf{g}_0)$-module with highest weight $\lambda$. 

\begin{thm}\label{decomposition as a classical crystal} For $(k,l)\in H^{\pmb{\diamondsuit}}$, we have the following decomposition of $C({\bf v}_{k,l})$ as a $U_q(\mf{g}_0)$-crystal :
\begin{itemize}
\item[(1)] If $k=0$, then $C({\bf v}_{k,l})\cong B_0(0)$.

\item[(2)] If $1\leq k\leq n$, then
\begin{equation*}
\begin{split}
C({\bf v}_{k,l})\cong  &
\begin{cases}
 B_0({\rm cl}(\varpi_k)), & \text{for \ $\pmb{\diamondsuit}=(\, {\def\lr#1{\multicolumn{1}{|@{\hspace{.6ex}}c@{\hspace{.6ex}}|}{\raisebox{-.08ex}{$#1$}}}\raisebox{-.1ex}
{$\scalebox{0.35}{\begin{array}[b]{cc}
\cline{1-1}\cline{2-2}
\lr{\ \ \ }&\lr{\ \ \ }\\ 
\cline{1-1}\cline{2-2}
\end{array}}$}}\, ,\, 
{\def\lr#1{\multicolumn{1}{|@{\hspace{.6ex}}c@{\hspace{.6ex}}|}{\raisebox{-.08ex}{$#1$}}}\raisebox{-.1ex}
{$\scalebox{0.35}{\begin{array}[b]{cc}
\cline{1-1}\cline{2-2}
\lr{\ \ \ }&\lr{\ \ \ }\\ 
\cline{1-1}\cline{2-2}
\end{array}}$}}\,)$},\\
\bigsqcup_{i=0}^k B_0({\rm cl}(\varpi_{k-i})), & \text{for \ $\pmb{\diamondsuit}=(\, {\def\lr#1{\multicolumn{1}{|@{\hspace{.6ex}}c@{\hspace{.6ex}}|}{\raisebox{-.3ex}{$#1$}}}\raisebox{-.1ex}
{$\scalebox{0.35}{\begin{array}[b]{c}
\cline{1-1}
\lr{\ \ \ }\\ 
\cline{1-1}
\end{array}}$}}\, ,\, {\def\lr#1{\multicolumn{1}{|@{\hspace{.6ex}}c@{\hspace{.6ex}}|}{\raisebox{-.08ex}{$#1$}}}\raisebox{-.1ex}
{$\scalebox{0.35}{\begin{array}[b]{cc}
\cline{1-1}\cline{2-2}
\lr{\ \ \ }&\lr{\ \ \ }\\ 
\cline{1-1}\cline{2-2}
\end{array}}$}} \, )$},\\
B_0({\rm cl}(\varpi_k)), & \text{for \ $\pmb{\diamondsuit}=(\, {\def\lr#1{\multicolumn{1}{|@{\hspace{.6ex}}c@{\hspace{.6ex}}|}{\raisebox{-.08ex}{$#1$}}}\raisebox{-.1ex}
{$\scalebox{0.35}{\begin{array}[b]{cc}
\cline{1-1}\cline{2-2}
\lr{\ \ \ }&\lr{\ \ \ }\\ 
\cline{1-1}\cline{2-2}
\end{array}}$}} \,,\,{\def\lr#1{\multicolumn{1}{|@{\hspace{.6ex}}c@{\hspace{.6ex}}|}{\raisebox{-.3ex}{$#1$}}}\raisebox{-.1ex}
{$\scalebox{0.35}{\begin{array}[b]{c}
\cline{1-1}
\lr{\ \ \ }\\ 
\cline{1-1}
\end{array}}$}}\, )$},\\
\bigsqcup_{i=0}^{[\frac{k}{2}]} B_0({\rm cl}(\varpi_{k-2i}))^{\oplus 2} , & \text{for \ $\pmb{\diamondsuit}=(\, {\def\lr#1{\multicolumn{1}{|@{\hspace{.6ex}}c@{\hspace{.6ex}}|}{\raisebox{-.3ex}{$#1$}}}\raisebox{-.8ex}
{$\scalebox{0.35}{\begin{array}[b]{c}
\cline{1-1}
\lr{\ \ \ }\\ 
\cline{1-1}
\lr{\ \ \ }\\
\cline{1-1}
\end{array}}$}}\, ,\, {\def\lr#1{\multicolumn{1}{|@{\hspace{.6ex}}c@{\hspace{.6ex}}|}{\raisebox{-.08ex}{$#1$}}}\raisebox{-.1ex}
{$\scalebox{0.35}{\begin{array}[b]{cc}
\cline{1-1}\cline{2-2}
\lr{\ \ \ }&\lr{\ \ \ }\\ 
\cline{1-1}\cline{2-2}
\end{array}}$}} \,  )$ with $k\neq n$},\\
\bigsqcup_{i=0}^{[\frac{n}{2}]} B_0({\rm cl}(\varpi_{n-2i})), & \text{for \ $\pmb{\diamondsuit}=(\, {\def\lr#1{\multicolumn{1}{|@{\hspace{.6ex}}c@{\hspace{.6ex}}|}{\raisebox{-.3ex}{$#1$}}}\raisebox{-.8ex}
{$\scalebox{0.35}{\begin{array}[b]{c}
\cline{1-1}
\lr{\ \ \ }\\ 
\cline{1-1}
\lr{\ \ \ }\\
\cline{1-1}
\end{array}}$}}\, ,\, {\def\lr#1{\multicolumn{1}{|@{\hspace{.6ex}}c@{\hspace{.6ex}}|}{\raisebox{-.08ex}{$#1$}}}\raisebox{-.1ex}
{$\scalebox{0.35}{\begin{array}[b]{cc}
\cline{1-1}\cline{2-2}
\lr{\ \ \ }&\lr{\ \ \ }\\ 
\cline{1-1}\cline{2-2}
\end{array}}$}} \,  )$ with $k=n$},\\
\end{cases}
\end{split}
\end{equation*}
where $B^{\oplus 2}=B \sqcup B$ for a crystal $B$ and $\varpi_0=0$.
\end{itemize}
\end{thm}
\pf (1) It is clear since $C({\bf v}_{0,l})=\{{\bf v}_{0,l}\}$ and ${\rm wt}({\bf v}_{0,l})=0$ by \eqref{weight of v_{k,l}}. 

(2) \textsc{Case 1}. Suppose that $\pmb{\diamondsuit}=(\, {\def\lr#1{\multicolumn{1}{|@{\hspace{.4ex}}c@{\hspace{.4ex}}|}{\raisebox{-.3ex}{$#1$}}}\raisebox{-.2ex}
{$\scalebox{0.5}{\begin{array}[b]{cc}
\cline{1-1}\cline{2-2}
\lr{\ \ \ }&\lr{\ \ \ }\\ 
\cline{1-1}\cline{2-2}
\end{array}}$}}\, ,\, 
{\def\lr#1{\multicolumn{1}{|@{\hspace{.4ex}}c@{\hspace{.4ex}}|}{\raisebox{-.3ex}{$#1$}}}\raisebox{-.2ex}
{$\scalebox{0.5}{\begin{array}[b]{cc}
\cline{1-1}\cline{2-2}
\lr{\ \ \ }&\lr{\ \ \ }\\ 
\cline{1-1}\cline{2-2}
\end{array}}$}}\,)$.
Let ${\bf m}\in C({\bf v}_{k,l})$ be given such that $\te_i{\bf m}=0$ for $i\in I_0$. By the same argument as in Proposition \ref{decomposition of M}, ${\bf m}={\bf v}_{k',l'}$ for some $k'$ and $l'$, which implies that ${\bf m}\in C({\bf v}_{k',l'})$.  Now from the decomposition of ${\bf M}$ as a $U_q'(\mf{g})$-crystal in Proposition \ref{decomposition of M}, it follows that $k'=k$ and $l'=l$. Therefore, $C({\bf v}_{k,l})$ is the connected as a $U_q({\mf g}_0)$-crystal. Since ${\bf M}$ is a regular crystal, $C({\bf v}_{k,l})$ is isomorphic to $B_0({\rm cl}(\varpi_k))$ as a $U_q({\mf g}_0)$-crystal by \eqref{weight of v_{k,l}}. The proof for $\pmb{\diamondsuit}=(\, {\def\lr#1{\multicolumn{1}{|@{\hspace{.4ex}}c@{\hspace{.4ex}}|}{\raisebox{-.3ex}{$#1$}}}\raisebox{-.2ex}
{$\scalebox{0.5}{\begin{array}[b]{cc}
\cline{1-1}\cline{2-2}
\lr{\ \ \ }&\lr{\ \ \ }\\ 
\cline{1-1}\cline{2-2}
\end{array}}$}}\,,\,{\def\lr#1{\multicolumn{1}{|@{\hspace{.4ex}}c@{\hspace{.4ex}}|}{\raisebox{-.3ex}{$#1$}}}\raisebox{-.2ex}
{$\scalebox{0.5}{\begin{array}[b]{c}
\cline{1-1}
\lr{\ \ \ }\\ 
\cline{1-1}
\end{array}}$}}\, )$ is almost the same.

\textsc{Case 2}. Suppose that $\pmb{\diamondsuit}=(\, {\def\lr#1{\multicolumn{1}{|@{\hspace{.4ex}}c@{\hspace{.4ex}}|}{\raisebox{-.3ex}{$#1$}}}\raisebox{-.2ex}
{$\scalebox{0.5}{\begin{array}[b]{c}
\cline{1-1}
\lr{\ \ \ }\\ 
\cline{1-1}
\end{array}}$}}\, ,\, {\def\lr#1{\multicolumn{1}{|@{\hspace{.4ex}}c@{\hspace{.4ex}}|}{\raisebox{-.3ex}{$#1$}}}\raisebox{-.2ex}
{$\scalebox{0.5}{\begin{array}[b]{cc}
\cline{1-1}\cline{2-2}
\lr{\ \ \ }&\lr{\ \ \ }\\ 
\cline{1-1}\cline{2-2}
\end{array}}$}} \, )$. Let ${\bf m}\in C({\bf v}_{k,n-k})$ be given such that $\te_i{\bf m}=0$ for $i\in I_0$.
As in Case 1, we have ${\bf m}={\bf v}_{k',l'}$ for some $k'$ and $l'$. We see that   ${\bf m}\in C({\bf v}_{n-l', l'})$  by applying $\tf_i$'s to ${\bf m}$ for $i\in \{  n-l'-1,\ldots, 1, 0\}$, and then $l'=n-k$ by Proposition \ref{decomposition of M}. Hence ${\bf m}={\bf v}_{k',n-k}$ with $0\leq k'\leq k$. Conversely, for $0\leq k'\leq k$, we have ${\bf v}_{k',n-k}\in C({\bf v}_{k,n-k})$  by applying $\tf_i$'s to ${\bf v}_{k',n-k}$ for $i\in \{  k',\ldots, 1, 0\}$. Hence 
\begin{equation}\label{classical decomposition-1}
C({\bf v}_{k,n-k})=\bigsqcup_{k'=0}^kC_0({\bf v}_{k',n-k}),
\end{equation}
where $C_0({\bf m})$ denotes the connected component of ${\bf m}$ as a $U_q({\mf g}_0)$-crystal. Finally, $C_0({\bf v}_{k',n-k})$ is isomorphic to $B_0({\rm cl}(\varpi_{k'}))$ by \eqref{weight of v_{k,l}}. This proves the decomposition of $C({\bf v}_{k,n-k})$. 

\textsc{Case 3}. Suppose that $\pmb{\diamondsuit}=(\, {\def\lr#1{\multicolumn{1}{|@{\hspace{.6ex}}c@{\hspace{.6ex}}|}{\raisebox{-.3ex}{$#1$}}}\raisebox{-.7ex}
{$\scalebox{0.45}{\begin{array}[b]{c}
\cline{1-1}
\lr{\ \ \ }\\ 
\cline{1-1}
\lr{\ \ \ }\\
\cline{1-1}
\end{array}}$}}\, ,\, {\def\lr#1{\multicolumn{1}{|@{\hspace{.4ex}}c@{\hspace{.4ex}}|}{\raisebox{-.3ex}{$#1$}}}\raisebox{-.2ex}
{$\scalebox{0.5}{\begin{array}[b]{cc}
\cline{1-1}\cline{2-2}
\lr{\ \ \ }&\lr{\ \ \ }\\ 
\cline{1-1}\cline{2-2}
\end{array}}$}} \,  )$. Similarly, from the argument in the proof of Proposition \ref{decomposition of M}, we see that if $k\neq n$, then
\begin{equation}\label{2 fold decomposition of C}
C({\bf v}_{k,n-k})=\bigsqcup_{i=0}^{[\frac{k}{2}]}\left( C_0({\bf v}_{k-2i ,n-k}) \sqcup C_0({\bf v}_{k-2i ,n-k-1}) \right),
\end{equation}
where both $C_0({\bf v}_{k-2i ,n-k})$ are $C_0({\bf v}_{k-2i ,n-k-1})$ are isomorphic to $B_0({\rm cl}(\varpi_{k-2i}))$. Also, if $k=n$, then we have 
\begin{equation}\label{classical decomposition-2}
C({\bf v}_{n,0})=\bigsqcup_{i=0}^{[\frac{n}{2}]} C_0({\bf v}_{n-2i ,0}),
\end{equation}
where $C_0({\bf v}_{n-2i ,0})\cong B_0({\rm cl}(\varpi_{n-2i}))$. The proof completes.
\qed
\vskip 3mm

\subsection{Order 2 symmetry of $A_{2n-1}^{(2)}$-crystals}
Let us introduce a symmetry of the crystal $C({\bf v}_{k,n-k})$ of type $\pmb{\diamondsuit}=(\, {\def\lr#1{\multicolumn{1}{|@{\hspace{.6ex}}c@{\hspace{.6ex}}|}{\raisebox{-.3ex}{$#1$}}}\raisebox{-.7ex}
{$\scalebox{0.45}{\begin{array}[b]{c}
\cline{1-1}
\lr{\ \ \ }\\ 
\cline{1-1}
\lr{\ \ \ }\\
\cline{1-1}
\end{array}}$}}\, ,\, {\def\lr#1{\multicolumn{1}{|@{\hspace{.4ex}}c@{\hspace{.4ex}}|}{\raisebox{-.3ex}{$#1$}}}\raisebox{-.2ex}
{$\scalebox{0.5}{\begin{array}[b]{cc}
\cline{1-1}\cline{2-2}
\lr{\ \ \ }&\lr{\ \ \ }\\ 
\cline{1-1}\cline{2-2}
\end{array}}$}} \,  )$ with $k\neq n$, which will play an important rule in the next section.

\begin{lem}\label{sign decomposition}
For ${\bf m}\in C({\bf v}_{k,n-k})$, we have $\sigma({\bf m})=(2i,n-k)$ or $(2i+1,n-k-1)$ for some $0\leq i\leq [\frac{k}{2}]$. 
\end{lem}
\pf Note that $\sigma({\bf v}_{k-2i ,n-k})=(2i,n-k)$ and $\sigma({\bf v}_{k-2i ,n-k-1})=(2i+1,n-k-1)$ for $0\leq i\leq [\frac{k}{2}]$. Note that $\sigma$ is constant on each connected component in $C({\bf v}_{k,n-k})$ as a $U_q(\mf{g}_0)$-crystal (see Case 1 in the proof of Proposition \ref{decomposition of M}). Hence the claim follows from \eqref{2 fold decomposition of C}.
\qed\vskip 2mm

By Lemma \ref{sign decomposition}, we have $C({\bf v}_{k,n-k})=C({\bf v}_{k,n-k})^+\sqcup C({\bf v}_{k,n-k})^-$, where
\begin{equation*}
\begin{split}
C({\bf v}_{k,n-k})^+&=\{\,{\bf m}\in C({\bf v}_{k,n-k})\,|\,\varphi({\bf m})=n-k \,\}=\bigsqcup_{i=0}^{[\frac{k}{2}]}  C_0({\bf v}_{k-2i ,n-k}),\\
C({\bf v}_{k,n-k})^-&=\{\,{\bf m}\in C({\bf v}_{k,n-k})\,|\,\varphi({\bf m})=n-k-1 \,\}=\bigsqcup_{i=0}^{[\frac{k}{2}]}  C_0({\bf v}_{k-2i ,n-k-1}).
\end{split}
\end{equation*}
 Now, define a map 
\begin{equation*}
\varsigma : C({\bf v}_{k,n-k}) \longrightarrow C({\bf v}_{k,n-k})
\end{equation*}
by $\varsigma({\bf m})=\td{F}{\bf m}$ (resp. $\td{E}{\bf m}$) if ${\bf m}\in C({\bf v}_{k,n-k})^+$ (resp. $C({\bf v}_{k,n-k})^-$). By definition, $\varsigma^2({\bf m})={\bf m}$ for ${\bf m}\in C({\bf v}_{k,n-k})$.

\begin{ex}{\rm
Let
$${\bf m}=
\begin{bmatrix}
\,1 & 0 \,\\
\,1 & 1 \,\\
\,0 & 1 \,
\end{bmatrix}\in C({\bf v}_{1,1}).$$
Then $\sigma({\bf m})=(1,1)$ and ${\bf m}\in C({\bf v}_{1,1})^-$. Hence 
$$\varsigma({\bf m})=\td{E}{\bf m}=
\begin{bmatrix}
\,1 & 0 \,\\
\,1 & 1 \,\\
\,1 & 0 \,
\end{bmatrix}\in C({\bf v}_{1,1})^+.$$

}
\end{ex}

\begin{prop}\label{varsigma}
$\varsigma$ commutes with $\te_i$ and $\tf_i$ for $i\in I$. Hence $\varsigma$ is an isomorphism of $U'_q({\mf g})$-crystals.
\end{prop}
\pf It is straightforward to check that $\varsigma$ commutes with $\te_0$ and $\tf_0$. Since $\tilde{E}$ and $\tilde{F}$ (and hence $\varsigma$) commute with $\te_i$ and $\tf_i$ for $i\in I_0$, this proves the claim.  \qed\vskip 2mm

We put $C({\bf v}_{k,n-k})/\langle \varsigma \rangle=\{\,{\bf m}+\varsigma({\bf m})\,|\,{\bf m}\in C({\bf v}_{k,n-k})\,\}$, which has a well-defined $I$-oriented graph structure  induced from $C({\bf v}_{k,n-k})$ (see Figure \ref{Graph C}).

\section{Realization of level zero fundamental representations}
\subsection{$W(\varpi_{n})$ or $W(\varpi_{n-1})$ of type $B_n^{(1)}$, $D_{n}^{(1)}$, $D_{n+1}^{(2)}$}
Suppose that ${\mf g}$ is of type $\pmb{\diamondsuit}=(\diamondsuit_0, \diamondsuit_n)$ with 
$\diamondsuit_0$, $\diamondsuit_n \neq{\def\lr#1{\multicolumn{1}{|@{\hspace{.4ex}}c@{\hspace{.4ex}}|}{\raisebox{-.3ex}{$#1$}}}\raisebox{-.2ex}
{$\scalebox{0.5}{\begin{array}[b]{cc}
\cline{1-1}\cline{2-2}
\lr{\ \ \ }&\lr{\ \ \ }\\ 
\cline{1-1}\cline{2-2}
\end{array}}$}}$\ , and ${\bf v}_{k}$ is a $U_q({\mf g}_0)$-highest weight element in ${\bf B}$ ($k=n, n-1$) (see  \eqref{decomposition of B}). As a $U_q({\mf g}_0)$-crystal, $C({\bf v}_k)$ is isomorphic to $B_0({\rm cl}(\varpi_k))$, the crystal of a spin representation.
Let $v_n=|0\rangle$ and $v_{n-1}=\psi_{\ov{n}}|0\rangle$. Then ${v_k}$ is a $U_q(\mf{g}_0)$-highest weight vector in $L(\E)$  such that $v_{k}\equiv {\bf v}_{k}$ or $\psi_{{\bf v}_k}|0\rangle  \pmod{q_s L(\E)}$. Let 
\begin{equation}
W_k=U'_q({\mf g}) v_k.
\end{equation}
By \eqref{decomposition of B}, $( {L}(\E)\cap W_k,C({\bf v}_k))$ is a crystal base of $W_k$ and 
\begin{equation}\label{decomposition of E}
\E=
\begin{cases}
W_n \oplus W_{n-1}, & \text{if \ $\pmb{\diamondsuit}= (\, {\def\lr#1{\multicolumn{1}{|@{\hspace{.6ex}}c@{\hspace{.6ex}}|}{\raisebox{-.3ex}{$#1$}}}\raisebox{-.8ex}
{$\scalebox{0.35}{\begin{array}[b]{c}
\cline{1-1}
\lr{\ \ \ }\\ 
\cline{1-1}
\lr{\ \ \ }\\
\cline{1-1}
\end{array}}$}}\, ,\, {\def\lr#1{\multicolumn{1}{|@{\hspace{.6ex}}c@{\hspace{.6ex}}|}{\raisebox{-.3ex}{$#1$}}}\raisebox{-.8ex}
{$\scalebox{0.35}{\begin{array}[b]{c}
\cline{1-1}
\lr{\ \ \ }\\ 
\cline{1-1}
\lr{\ \ \ }\\
\cline{1-1}
\end{array}}$}} \, )$},\\
W_n, & \text{otherwise}.\\
\end{cases}\end{equation}
Since all the weights of $W_k$ are minuscule and hence contained in the convex hull of $W{\rm cl}(\varpi_k)$,  $v_k$ is an extremal weight vector by \cite[Theorem 5.3]{Kas02}. Since ${\rm dim}(W_k)_{{\rm cl}(\varpi_k)}=1$ and $C({\bf v}_k)$ is connected, $W_k$ is irreducible (cf.\cite[Lemma 2.7]{BKK}) and therefore $W_k\cong W(\varpi_k)_{a_k}$ 
for some $a_k\in K\setminus\{0\}$ by Theorem \ref{Fundamental Repn} (6). Choose $w\in W$
such that $w(\varpi_k)=\varpi_k+\delta$. Note that if $w=s_{i_1}\ldots s_{i_t}$, then $S_w v_k= x_{i_1}^{(m_1)}\ldots x_{i_t}^{(m_t)}v_k$, where $x\in \{e,f\}$, and $m_s\in\{0,1\}$.  Then by Propositions \ref{rho-1} and \ref{local-1}, it is not difficult  to see that $S_w v_k = v_k$. Since $S_w u_{\varpi_k} = z_k u_{\varpi_k}$ in $V(\varpi_k)$ and $W(\varpi_k)_a=V(\varpi_k)/(z_k-a)V(\varpi_k)$, we have $a_k=1$. Therefore
\begin{equation}\label{W_k}
W_k\cong W(\varpi_k).
\end{equation}
We remark that the construction of $W(\varpi_k)$ ($k=n-1,n$) is  already well-known. The decomposition of $\E^{\otimes 2}$ follows from \eqref{decomposition of E}, \eqref{W_k}, and simplicity of tensor product of $W(\varpi_k)$'s \cite[Theorem 9.2]{Kas02}.

\subsection{$W(\varpi_k)$ ($1\leq k\leq n$) of type $C_{n}^{(1)}$, $A_{2n}^{(2)}$, $A_{2n}^{(2)\dagger}$, $A_{2n-1}^{(2)}$} Suppose that ${\mf g}$ is of type $\pmb{\diamondsuit}=(\diamondsuit_0, \diamondsuit_n)$ with 
  $\diamondsuit_0$ or $\diamondsuit_n ={\def\lr#1{\multicolumn{1}{|@{\hspace{.4ex}}c@{\hspace{.4ex}}|}{\raisebox{-.3ex}{$#1$}}}\raisebox{-.2ex}
{$\scalebox{0.5}{\begin{array}[b]{cc}
\cline{1-1}\cline{2-2}
\lr{\ \ \ }&\lr{\ \ \ }\\ 
\cline{1-1}\cline{2-2}
\end{array}}$}}$\ .
Let  $(k,l)\in H^{\pmb{\diamondsuit}}$ be given.
Since ${\bf v}_{k,l}\in {\bf M}$ is a $U_q({\mf g}_0)$-highest weight crystal element, there exists a  $U_q(\mf{g}_0)$-highest weight vector  $v_{k,l}$ in $L(\E^{\otimes 2})$ such that $v_{k,l}\equiv {\bf v}_{k,l}$ or $\psi_{{\bf v}_{k,l}^{(1)}}|0\rangle\otimes \psi_{{\bf v}_{k,l}^{(2)}}|0\rangle \pmod{q {L}(\E^{\otimes 2})}$. We may also assume that the coefficient of $\psi_{{\bf v}_{k,l}^{(1)}}|0\rangle\otimes \psi_{{\bf v}_{k,l}^{(2)}}|0\rangle$ in $v_{k,l}$ is $1$ after multiplication by an invertible element in $\mathbb{A}$. 
We put
\begin{equation}
W_{k,l}=U'_q({\mf g})v_{k,l}.
\end{equation}
For 
 $\pmb{\diamondsuit}=(\, {\def\lr#1{\multicolumn{1}{|@{\hspace{.6ex}}c@{\hspace{.6ex}}|}{\raisebox{-.3ex}{$#1$}}}\raisebox{-.7ex}
{$\scalebox{0.45}{\begin{array}[b]{c}
\cline{1-1}
\lr{\ \ \ }\\ 
\cline{1-1}
\lr{\ \ \ }\\
\cline{1-1}
\end{array}}$}}\, ,\, {\def\lr#1{\multicolumn{1}{|@{\hspace{.4ex}}c@{\hspace{.4ex}}|}{\raisebox{-.3ex}{$#1$}}}\raisebox{-.2ex}
{$\scalebox{0.5}{\begin{array}[b]{cc}
\cline{1-1}\cline{2-2}
\lr{\ \ \ }&\lr{\ \ \ }\\ 
\cline{1-1}\cline{2-2}
\end{array}}$}} \,  )$ with $k\not\in\{ 0, n\}$,
put $v_{k,n-k}^\pm=v_{k,n-k}\pm v_{k,n-k-1}$ and
\begin{equation}
W_{k,n-k}^\pm = U'_q(\mf{g}) v_{k,n-k}^\pm.
\end{equation}
Here $v_{k,n-k-1}$ is also assumed to be a $U_q(\mf{g}_0)$-highest weight vector such that $v_{k,n-k-1}\in \psi_{{\bf v}_{k,n-k-1}^{(1)}}|0\rangle\otimes \psi_{{\bf v}_{k,n-k-1}^{(2)}}|0\rangle +q {L}(\E^{\otimes 2})$.

\begin{prop}\label{decomposition of E^2} 
We have $$\E^{\otimes 2}=\bigoplus_{(k,l)\in  H^{\pmb{\diamondsuit}}}W_{k,l},$$
and $W_{k,l}$ has a polarizable crystal base $(L(W_{k,l}),B(W_{k,l}))$ with $L(W_{k,l})=L(\E^{\otimes 2})\cap W_{k,l}$ and $B(W_{k,l})=C({\bf v}_{k,l})$.
\end{prop}
\pf Note that $L(W_{k,l})$ is invariant under $\te_i$ and $\tf_i$ for $i\in I$, and $C({\bf v}_{k,l})$ is linearly independent subset of $L(W_{k,l})/q_sL(W_{k,l})$. By Proposition \ref{decomposition of M}, we have $L(\E^{\otimes 2})=\bigoplus_{(k,l)\in H^{\pmb{\diamondsuit}}}L(W_{k,l})$, which implies that $\E^{\otimes 2}=\bigoplus_{(k,l)\in  H^{\pmb{\diamondsuit}}}W_{k,l}$, and $(L(W_{k,l}), C({\bf v}_{k,l}))$ is a crystal base of $W_{k,l}$. The polarizability follows from that of $\E^{\otimes 2}$.
\qed

\begin{lem}\label{extremality of v_{k,l}}
For $(k,l)\in H^{\pmb{\diamondsuit}}$, $v_{k,l}$ is an extremal weight vector of weight ${\rm cl}(\varpi_k)$.
\end{lem}
\pf  By  Proposition \ref{decomposition of E^2}, we have  $B(W_{k,l})=C({\bf v}_{k,l})$. 
By Theorem \ref{decomposition as a classical crystal}, we see that  the weights of $W_{k,l}$ or $C({\bf v}_{k,l})$ belong to the convex hull of $W{\rm cl}(\varpi_k)$. Then by \cite[Theorem 5.3]{Kas02}, $v_{k,l}$ is an extremal weight vector.
\qed

\begin{lem}\label{S_w v_{k,l}}
For $w\in W$ and $(k,l)\in H^{\pmb{\diamondsuit}}$, put {\rm ${\bf m}={\texttt S}_w{\bf v}_{k,l}$}. Then {\rm
\begin{equation*}
S_w v_{k,l} =\psp_{{\bf m}^{(1)}}|0\rangle \otimes \psp_{{\bf m}^{(2)}}|0\rangle +\sum_{\substack{{\bf m}'\in {\bf M}\setminus\{{\texttt S}_w{\bf v}_{k,l}\}\\ {\rm wt}({\bf m}')=w ({\rm cl}(\varpi_k)) }}a_{{\bf m}'}\psp_{{\bf m}^{'(1)}}|0\rangle \otimes \psp_{{\bf m}^{'(2)}}|0\rangle,
\end{equation*} }
for some $a_{{\bf m}'}\in q_s \mathbb{A}$.
\end{lem}
\pf When $w=1$, $v_{k,l}$ satisfies the above property by our choice of $v_{k,l}$. In general, it follows directly from the descriptions of $e_i$ and $f_i$ ($i\in I$) on $\E^{\otimes 2}$ in  Propositions \ref{rho-1}, \ref{local-1}, and \ref{local-2}, and the induction on the length of $w$ (note that $S_w v_{k,l}$ is an extremal weight vector for all $w\in W$).
\qed

\begin{thm}\label{main theorem}
Let $(k,l)\in H^{\pmb{\diamondsuit}}$ be given.
\begin{itemize}
\item[(1)] If $\pmb{\diamondsuit}\neq (\, {\def\lr#1{\multicolumn{1}{|@{\hspace{.4ex}}c@{\hspace{.4ex}}|}{\raisebox{-.3ex}{$#1$}}}\raisebox{-.7ex}
{$\scalebox{0.45}{\begin{array}[b]{c}
\cline{1-1}
\lr{\ \ \ }\\ 
\cline{1-1}
\lr{\ \ \ }\\
\cline{1-1}
\end{array}}$}}\, ,\, {\def\lr#1{\multicolumn{1}{|@{\hspace{.3ex}}c@{\hspace{.3ex}}|}{\raisebox{-.3ex}{$#1$}}}\raisebox{-.2ex}
{$\scalebox{0.5}{\begin{array}[b]{cc}
\cline{1-1}\cline{2-2}
\lr{\ \ \ }&\lr{\ \ \ }\\ 
\cline{1-1}\cline{2-2}
\end{array}}$}} \,  )$ or
$\pmb{\diamondsuit}=(\, {\def\lr#1{\multicolumn{1}{|@{\hspace{.4ex}}c@{\hspace{.4ex}}|}{\raisebox{-.3ex}{$#1$}}}\raisebox{-.7ex}
{$\scalebox{0.45}{\begin{array}[b]{c}
\cline{1-1}
\lr{\ \ \ }\\ 
\cline{1-1}
\lr{\ \ \ }\\
\cline{1-1}
\end{array}}$}}\, ,\, {\def\lr#1{\multicolumn{1}{|@{\hspace{.4ex}}c@{\hspace{.4ex}}|}{\raisebox{-.3ex}{$#1$}}}\raisebox{-.2ex}
{$\scalebox{0.5}{\begin{array}[b]{cc}
\cline{1-1}\cline{2-2}
\lr{\ \ \ }&\lr{\ \ \ }\\ 
\cline{1-1}\cline{2-2}
\end{array}}$}} \,  )$ with $k\in\{ 0, n\}$, then
$$W_{k,l}\cong W(\varpi_k).$$

\item[(2)] If  $\pmb{\diamondsuit}=(\, {\def\lr#1{\multicolumn{1}{|@{\hspace{.4ex}}c@{\hspace{.4ex}}|}{\raisebox{-.3ex}{$#1$}}}\raisebox{-.7ex}
{$\scalebox{0.45}{\begin{array}[b]{c}
\cline{1-1}
\lr{\ \ \ }\\ 
\cline{1-1}
\lr{\ \ \ }\\
\cline{1-1}
\end{array}}$}}\, ,\, {\def\lr#1{\multicolumn{1}{|@{\hspace{.4ex}}c@{\hspace{.4ex}}|}{\raisebox{-.3ex}{$#1$}}}\raisebox{-.2ex}
{$\scalebox{0.5}{\begin{array}[b]{cc}
\cline{1-1}\cline{2-2}
\lr{\ \ \ }&\lr{\ \ \ }\\ 
\cline{1-1}\cline{2-2}
\end{array}}$}} \,  )$ and $k\not\in\{ 0, n\}$, then  $W_{k,n-k}=W_{k,n-k}^+\oplus W_{k,n-k}^-$ and
$$W_{k,n-k}^\pm\cong W(\varpi_k)_{\pm 1}$$
\end{itemize}
Here, $W(\varpi_0)=W(0)$ is the trivial $U'_q(\mf{g})$-module of one dimension.  
\end{thm}
\pf (1) Suppose that either $\pmb{\diamondsuit}\neq (\, {\def\lr#1{\multicolumn{1}{|@{\hspace{.6ex}}c@{\hspace{.6ex}}|}{\raisebox{-.3ex}{$#1$}}}\raisebox{-.7ex}
{$\scalebox{0.45}{\begin{array}[b]{c}
\cline{1-1}
\lr{\ \ \ }\\ 
\cline{1-1}
\lr{\ \ \ }\\
\cline{1-1}
\end{array}}$}}\, ,\, {\def\lr#1{\multicolumn{1}{|@{\hspace{.4ex}}c@{\hspace{.4ex}}|}{\raisebox{-.3ex}{$#1$}}}\raisebox{-.2ex}
{$\scalebox{0.5}{\begin{array}[b]{cc}
\cline{1-1}\cline{2-2}
\lr{\ \ \ }&\lr{\ \ \ }\\ 
\cline{1-1}\cline{2-2}
\end{array}}$}} \,  )$ or
$\pmb{\diamondsuit}=(\, {\def\lr#1{\multicolumn{1}{|@{\hspace{.6ex}}c@{\hspace{.6ex}}|}{\raisebox{-.3ex}{$#1$}}}\raisebox{-.7ex}
{$\scalebox{0.45}{\begin{array}[b]{c}
\cline{1-1}
\lr{\ \ \ }\\ 
\cline{1-1}
\lr{\ \ \ }\\
\cline{1-1}
\end{array}}$}}\, ,\, {\def\lr#1{\multicolumn{1}{|@{\hspace{.4ex}}c@{\hspace{.4ex}}|}{\raisebox{-.3ex}{$#1$}}}\raisebox{-.2ex}
{$\scalebox{0.5}{\begin{array}[b]{cc}
\cline{1-1}\cline{2-2}
\lr{\ \ \ }&\lr{\ \ \ }\\ 
\cline{1-1}\cline{2-2}
\end{array}}$}} \,  )$ with $k\in\{ 0, n\}$.
Recall that $B(W_{k,l})$ is connected and $\dim (W_{k,l})_{{\rm cl}(\varpi_k)}=\left\vert C({\bf v}_{k,l})_{{\rm cl}(\varpi_k)}\right\vert=1$ by Theorem \ref{decomposition as a classical crystal}. Hence $W_{k,l}$ is an irreducible $U'_q(\mf{g})$-module generated by $v_{k,l}$. When $k=0$, it is clear that $W_{k,l}$ the trivial $U'_q(\mf{g})$-module of one dimension. 
Assume that $k\neq 0$. By Lemma \ref{extremality of v_{k,l}} and Theorem \ref{Fundamental Repn} (6), $W_{k,l}\cong W(\varpi_k)_{a_{k,l}}$ for some $a_{k,l}\in K\setminus\{0\}$. Choose $w\in W$ such that $w(\varpi_k)=\varpi_k+d_k\delta$. Since $\dim  C({\bf v}_{k,l})_{{\rm cl}(\varpi_k)}=1$, we have ${\texttt S}_w {\bf v}_{k,l}={\bf v}_{k,l}$. By Lemma \ref{S_w v_{k,l}}, we have
\begin{equation*}
\begin{split}
S_w v_{k,l} &=\psi_{{\bf v}_{k,l}^{(1)}}|0\rangle\otimes \psi_{{\bf v}_{k,l}^{(2)}}|0\rangle +\sum_{\substack{{\bf m}\in {\bf M}\setminus\{{\bf v}_{k,l}\}\\ {\rm wt}({\bf m})={\rm cl}(\varpi_k)}}a_{\bf m} \psp_{{\bf m}^{(1)}}|0\rangle \otimes \psp_{{\bf m}^{(2)}}|0\rangle,
\end{split}
\end{equation*} 
for some $a_{\bf m}\in q_s\mathbb{A}$.
Since $S_w v_{k,l}=a_{k,l}v_{k,l}$, we have $a_{k,l}=1$, and hence $W_{k,l}$ is isomorphic to $W(\varpi_k)$.

(2) Suppose that 
 $\pmb{\diamondsuit}=(\, {\def\lr#1{\multicolumn{1}{|@{\hspace{.6ex}}c@{\hspace{.6ex}}|}{\raisebox{-.3ex}{$#1$}}}\raisebox{-.7ex}
{$\scalebox{0.45}{\begin{array}[b]{c}
\cline{1-1}
\lr{\ \ \ }\\ 
\cline{1-1}
\lr{\ \ \ }\\
\cline{1-1}
\end{array}}$}}\, ,\, {\def\lr#1{\multicolumn{1}{|@{\hspace{.4ex}}c@{\hspace{.4ex}}|}{\raisebox{-.3ex}{$#1$}}}\raisebox{-.2ex}
{$\scalebox{0.5}{\begin{array}[b]{cc}
\cline{1-1}\cline{2-2}
\lr{\ \ \ }&\lr{\ \ \ }\\ 
\cline{1-1}\cline{2-2}
\end{array}}$}} \,  )$ and $k\not\in \{0, n\}$.
Put
\begin{equation*}
\begin{split}
%W_{k,n-k}^\pm &= U'_q(\mf{g})(v_{k,n-k}\pm v_{k,n-k-1}),\\
L(W_{k,n-k}^\pm)&=L(W_{k,n-k})\cap W_{k,n-k}^\pm,\\
B(W_{k,n-k}^\pm)&=\{\,{\bf m}\pm\varsigma({\bf m})\,|\,{\bf m}\in C({\bf v}_{k,n-k})^+\,\}.
\end{split}
\end{equation*}
Note that $B(W_{k,n-k}^+)\sqcup B(W_{k,n-k}^-)$ is a $\mathbb{Q}$-basis of $L(W_{k,n-k})/qL(W_{k,n-k})$ since it is orthogonal. By Proposition \ref{varsigma}, $B(W_{k,n-k}^\pm)$ is a linearly independent subset of $L(W_{k,n-k}^\pm)/q_sL(W_{k,n-k}^\pm)$, which is invariant under $\te_i$ and $\tf_i$ for $i\in I$ up to multiplication by $\pm1$. This implies that $$W_{k,n-k}=W_{k,n-k}^+\oplus W_{k,n-k}^-,$$ and $(L(W_{k,n-k}^\pm),B(W_{k,n-k}^\pm))$ is a pseudo crystal base of $W_{k,n-k}^\pm$ (cf. \cite{KMN2}). Moreover, $B(W_{k,n-k}^\pm)\cong C({\bf v}_{k,n-k})/\langle \varsigma \rangle$ as an $I$-colored oriented graph.

Since $B(W_{k,n-k}^\pm)$ is connected with $\left\vert B(W_{k,n-k}^\pm)_{{\rm cl}(\varpi_k)}\right\vert=1$ and $v_{k,n-k}^\pm$ is an extremal weight vector, we have $W_{k,n-k}^\pm\cong W(\varpi_k)_{a^\pm_{k,n-k}}$ for some $a^\pm_{k,n-k}\in K\setminus\{0\}$ by the same arguments as in (1).

Finally, let $w\in W$ be given such that $w(\varpi_k)=\varpi_k+\delta$. Since ${\texttt S}_w ({\bf v}_{k,n-k}\pm {\bf v}_{k,n-k-1})=\pm ({\bf v}_{k,n-k}\pm {\bf v}_{k,n-k-1})$ (see \eqref{delta shifting w-2}), we have $S_w v_{k,n-k}^\pm =\pm v_{k,n-k}^\pm$ by Lemma \ref{S_w v_{k,l}} and hence  $a^\pm_{k,n-k}=\pm 1$. Therefore, $W_{k,n-k}^\pm$ is isomorphic to $W(\varpi_k)_{\pm 1}$. 
\qed\vskip 2mm

Let $B(W(\varpi_k))$ denote the crystal of $W(\varpi_k)$ for $k\in I_0$.

\begin{cor}\label{crystal realization}
For $(k,l)\in H^{\pmb{\diamondsuit}}$ with $k\neq 0$, we have
\begin{equation*}
B(W(\varpi_k))\cong
\begin{cases}
C({\bf v}_{k,n-k})/\langle \varsigma \rangle, & \text{if $\pmb{\diamondsuit}=(\, {\def\lr#1{\multicolumn{1}{|@{\hspace{.6ex}}c@{\hspace{.6ex}}|}{\raisebox{-.3ex}{$#1$}}}\raisebox{-.8ex}
{$\scalebox{0.35}{\begin{array}[b]{c}
\cline{1-1}
\lr{\ \ \ }\\ 
\cline{1-1}
\lr{\ \ \ }\\
\cline{1-1}
\end{array}}$}}\, ,\, {\def\lr#1{\multicolumn{1}{|@{\hspace{.6ex}}c@{\hspace{.6ex}}|}{\raisebox{-.08ex}{$#1$}}}\raisebox{-.1ex}
{$\scalebox{0.35}{\begin{array}[b]{cc}
\cline{1-1}\cline{2-2}
\lr{\ \ \ }&\lr{\ \ \ }\\ 
\cline{1-1}\cline{2-2}
\end{array}}$}} \, )$ and $k\neq n$}, \\
C({\bf v}_{k,l}), & \text{otherwise}.
\end{cases}
\end{equation*}
\end{cor}

\begin{cor}\label{decomposition of E^2-2} As a $U'_q(\mf{g})$-module, we have
\begin{equation*}
\E^{\otimes 2}\cong
\begin{cases}
\bigoplus_{0\leq k\leq n} W(\varpi_k)^{\oplus n-k+1}, & \text{if \ $\pmb{\diamondsuit}=(\, {\def\lr#1{\multicolumn{1}{|@{\hspace{.6ex}}c@{\hspace{.6ex}}|}{\raisebox{-.08ex}{$#1$}}}\raisebox{-.1ex}
{$\scalebox{0.35}{\begin{array}[b]{cc}
\cline{1-1}\cline{2-2}
\lr{\ \ \ }&\lr{\ \ \ }\\ 
\cline{1-1}\cline{2-2}
\end{array}}$}}\, ,\, 
{\def\lr#1{\multicolumn{1}{|@{\hspace{.6ex}}c@{\hspace{.6ex}}|}{\raisebox{-.08ex}{$#1$}}}\raisebox{-.1ex}
{$\scalebox{0.35}{\begin{array}[b]{cc}
\cline{1-1}\cline{2-2}
\lr{\ \ \ }&\lr{\ \ \ }\\ 
\cline{1-1}\cline{2-2}
\end{array}}$}}\,)$},\\
\bigoplus_{0\leq k\leq n} W(\varpi_k), & \text{if \ $\pmb{\diamondsuit}=(\, {\def\lr#1{\multicolumn{1}{|@{\hspace{.6ex}}c@{\hspace{.6ex}}|}{\raisebox{-.3ex}{$#1$}}}\raisebox{-.1ex}
{$\scalebox{0.35}{\begin{array}[b]{c}
\cline{1-1}
\lr{\ \ \ }\\ 
\cline{1-1}
\end{array}}$}}\, ,\, {\def\lr#1{\multicolumn{1}{|@{\hspace{.6ex}}c@{\hspace{.6ex}}|}{\raisebox{-.08ex}{$#1$}}}\raisebox{-.1ex}
{$\scalebox{0.35}{\begin{array}[b]{cc}
\cline{1-1}\cline{2-2}
\lr{\ \ \ }&\lr{\ \ \ }\\ 
\cline{1-1}\cline{2-2}
\end{array}}$}} \, )$ , $(\, {\def\lr#1{\multicolumn{1}{|@{\hspace{.6ex}}c@{\hspace{.6ex}}|}{\raisebox{-.08ex}{$#1$}}}\raisebox{-.1ex}
{$\scalebox{0.35}{\begin{array}[b]{cc}
\cline{1-1}\cline{2-2}
\lr{\ \ \ }&\lr{\ \ \ }\\ 
\cline{1-1}\cline{2-2}
\end{array}}$}} \,,\,{\def\lr#1{\multicolumn{1}{|@{\hspace{.6ex}}c@{\hspace{.6ex}}|}{\raisebox{-.3ex}{$#1$}}}\raisebox{-.1ex}
{$\scalebox{0.35}{\begin{array}[b]{c}
\cline{1-1}
\lr{\ \ \ }\\ 
\cline{1-1}
\end{array}}$}}\, )$},\\
W(\varpi_n)\oplus W(0)^{\oplus 2}\oplus \bigoplus_{\substack{1\leq k\leq n-1 \\ a=\pm1}}W(\varpi_k)_a, & \text{if \ $\pmb{\diamondsuit}=(\, {\def\lr#1{\multicolumn{1}{|@{\hspace{.6ex}}c@{\hspace{.6ex}}|}{\raisebox{-.3ex}{$#1$}}}\raisebox{-.8ex}
{$\scalebox{0.35}{\begin{array}[b]{c}
\cline{1-1}
\lr{\ \ \ }\\ 
\cline{1-1}
\lr{\ \ \ }\\
\cline{1-1}
\end{array}}$}}\, ,\, {\def\lr#1{\multicolumn{1}{|@{\hspace{.6ex}}c@{\hspace{.6ex}}|}{\raisebox{-.08ex}{$#1$}}}\raisebox{-.1ex}
{$\scalebox{0.35}{\begin{array}[b]{cc}
\cline{1-1}\cline{2-2}
\lr{\ \ \ }&\lr{\ \ \ }\\ 
\cline{1-1}\cline{2-2}
\end{array}}$}} \, )$.}
\end{cases}
\end{equation*}
\end{cor}

\begin{rem}{\rm 
By Theorem \ref{decomposition as a classical crystal}, we recover the decomposition of $B(W(\varpi_k))$ into $U_q(\mf{g}_0)$-crystals (see \cite{HN} for a more general case of ${\mf g}$ and $W(\varpi_k)$). For example, if $\pmb{\diamondsuit}=(\, {\def\lr#1{\multicolumn{1}{|@{\hspace{.6ex}}c@{\hspace{.6ex}}|}{\raisebox{-.3ex}{$#1$}}}\raisebox{-.8ex}
{$\scalebox{0.35}{\begin{array}[b]{c}
\cline{1-1}
\lr{\ \ \ }\\ 
\cline{1-1}
\lr{\ \ \ }\\
\cline{1-1}
\end{array}}$}}\, ,\, {\def\lr#1{\multicolumn{1}{|@{\hspace{.6ex}}c@{\hspace{.6ex}}|}{\raisebox{-.08ex}{$#1$}}}\raisebox{-.1ex}
{$\scalebox{0.35}{\begin{array}[b]{cc}
\cline{1-1}\cline{2-2}
\lr{\ \ \ }&\lr{\ \ \ }\\ 
\cline{1-1}\cline{2-2}
\end{array}}$}} \, )$ and $k\neq n$, then
\begin{equation*}
B(W(\varpi_k))\cong C({\bf v}_{k,n-k})/\langle \varsigma \rangle\cong \bigsqcup_{i=0}^{[\frac{k}{2}]} B_0({\rm cl}(\varpi_{k-2i})).
\end{equation*}  }
\end{rem}

\subsection{Description of $B(W(\varpi_k))$}
For $k\in \Z$, put $\langle k\rangle= \max\{k,0\}$.
By \cite[Proposition 2.1.1]{KN}, we have for ${\bf m}\in {\bf M}$,
\begin{equation*}
\begin{split}
\varepsilon({\bf m})&= \max\left\{\,\sum_{1\leq i\leq k}\langle m_{\ov{i}2}-m_{\ov{i}1}\rangle-\sum_{1\leq i <k} \langle m_{\ov{i}1}-m_{\ov{i}2} \rangle \,\Bigg|\,1\leq k\leq n\,\right\},\\
\varphi({\bf m})&=\max\left\{\,\sum_{k\leq i < n}\left\{\langle m_{\ov{i}1}-m_{\ov{i}2}\rangle - \langle m_{\ov{i+1}2}-m_{\ov{i+1}1}\rangle\right\}+\langle m_{\ov{n}1}-m_{\ov{n}2}\rangle\,\Bigg|\,1\leq k\leq n\,\right\}.
\end{split}
\end{equation*}
\vskip 2mm

By using the result \cite{K13-2} on the crystal of type $B_n$ and $C_n$, we have the following characterization of $B(W(\varpi_k))$ ($1\leq k\leq n$) of type $C_{n}^{(1)}$, $A_{2n}^{(2)}$, $A_{2n}^{(2)\dagger}$, $A_{2n-1}^{(2)}$ in terms of binary matrices ${\bf m}$ in ${\bf M}$ with constraints on $\sigma({\bf m})$. 

\begin{thm}\label{main theorem-2}
For $1\leq k\leq n$, we have the following.
\begin{itemize}
\item[(1)] If $\pmb{\diamondsuit}=(\, {\def\lr#1{\multicolumn{1}{|@{\hspace{.4ex}}c@{\hspace{.4ex}}|}{\raisebox{-.3ex}{$#1$}}}\raisebox{-.2ex}
{$\scalebox{0.5}{\begin{array}[b]{cc}
\cline{1-1}\cline{2-2}
\lr{\ \ \ }&\lr{\ \ \ }\\ 
\cline{1-1}\cline{2-2}
\end{array}}$}}\, ,\, 
{\def\lr#1{\multicolumn{1}{|@{\hspace{.4ex}}c@{\hspace{.4ex}}|}{\raisebox{-.3ex}{$#1$}}}\raisebox{-.2ex}
{$\scalebox{0.5}{\begin{array}[b]{cc}
\cline{1-1}\cline{2-2}
\lr{\ \ \ }&\lr{\ \ \ }\\ 
\cline{1-1}\cline{2-2}
\end{array}}$}}\,)$ or ${\mf g}=C_n^{(1)}$, then  
\begin{equation*}
B(W(\varpi_k))\cong\{\,{\bf m}\,|\,\sigma({\bf m})=(n-k,0)\,\}.
\end{equation*}

\item[(2)] If $\pmb{\diamondsuit}=(\, {\def\lr#1{\multicolumn{1}{|@{\hspace{.4ex}}c@{\hspace{.4ex}}|}{\raisebox{-.3ex}{$#1$}}}\raisebox{-.2ex}
{$\scalebox{0.5}{\begin{array}[b]{c}
\cline{1-1}
\lr{\ \ \ }\\ 
\cline{1-1}
\end{array}}$}}\, ,\, {\def\lr#1{\multicolumn{1}{|@{\hspace{.4ex}}c@{\hspace{.4ex}}|}{\raisebox{-.3ex}{$#1$}}}\raisebox{-.2ex}
{$\scalebox{0.5}{\begin{array}[b]{cc}
\cline{1-1}\cline{2-2}
\lr{\ \ \ }&\lr{\ \ \ }\\ 
\cline{1-1}\cline{2-2}
\end{array}}$}} \, )$ or ${\mf g}=A_{2n}^{(2)}$, then 
\begin{equation*}
B(W(\varpi_k))\cong\{\,{\bf m}\,|\,\sigma({\bf m})=(k-l,n-k) \text{ for $0\leq l\leq k$}\,\}.
\end{equation*}

\item[(3)] If $\pmb{\diamondsuit}=(\, {\def\lr#1{\multicolumn{1}{|@{\hspace{.4ex}}c@{\hspace{.4ex}}|}{\raisebox{-.3ex}{$#1$}}}\raisebox{-.2ex}
{$\scalebox{0.5}{\begin{array}[b]{cc}
\cline{1-1}\cline{2-2}
\lr{\ \ \ }&\lr{\ \ \ }\\ 
\cline{1-1}\cline{2-2}
\end{array}}$}}\,,\,{\def\lr#1{\multicolumn{1}{|@{\hspace{.4ex}}c@{\hspace{.4ex}}|}{\raisebox{-.3ex}{$#1$}}}\raisebox{-.2ex}
{$\scalebox{0.5}{\begin{array}[b]{c}
\cline{1-1}
\lr{\ \ \ }\\ 
\cline{1-1}
\end{array}}$}}\, )$ or ${\mf g}=A_{2n}^{(2)\dagger}$, then
\begin{equation*}
\begin{split}
& B(W(\varpi_k))\cong\left\{\,{\bf m}\, |\, \sigma({\bf m})=(n-k,l) \text{ for $0\leq l\leq k$}\right\}.
\end{split}
\end{equation*}

\item[(4)] If $\pmb{\diamondsuit}=(\, {\def\lr#1{\multicolumn{1}{|@{\hspace{.6ex}}c@{\hspace{.6ex}}|}{\raisebox{-.3ex}{$#1$}}}\raisebox{-.7ex}
{$\scalebox{0.45}{\begin{array}[b]{c}
\cline{1-1}
\lr{\ \ \ }\\ 
\cline{1-1}
\lr{\ \ \ }\\
\cline{1-1}
\end{array}}$}}\, ,\, {\def\lr#1{\multicolumn{1}{|@{\hspace{.4ex}}c@{\hspace{.4ex}}|}{\raisebox{-.3ex}{$#1$}}}\raisebox{-.2ex}
{$\scalebox{0.5}{\begin{array}[b]{cc}
\cline{1-1}\cline{2-2}
\lr{\ \ \ }&\lr{\ \ \ }\\ 
\cline{1-1}\cline{2-2}
\end{array}}$}} \,  )$ or ${\mf g}=A_{2n-1}^{(2)}$, then
\begin{equation*}
\begin{split}
B(W(\varpi_k))&\cong
\begin{cases}
\{\,{\bf m}\,|\,\sigma({\bf m})=(2l,0) \text{ for $0\leq l\leq [n/2]$}\,\},& \text{$(k=n)$},\\
\{\,{\bf m}+\varsigma({\bf m})\,|\,\sigma({\bf m})=(2l,n-k) \text{ for $0\leq l\leq [n/2]$}\,\}, & \text{$(k\neq n)$}.
\end{cases}
\end{split}
\end{equation*}

\end{itemize}
\end{thm}
\pf 
(1) We have $B(W(\varpi_k))\cong C({\bf v}_{k,0})= C_0({\bf v}_{k,0})$ by Corollary \ref{crystal realization} and Theorem \ref{decomposition as a classical crystal}. Recall that $\sigma({\bf m})=(n-k,0)=\sigma({\bf v}_{k,0})$ for ${\bf m}\in C({\bf v}_{k,0})$ since $\sigma$ is constant on $C_0({\bf v}_{k,0})$. 

Let $C=\{\,{\bf m}\in {\bf M}\,|\,\sigma({\bf m})=(n-k,0)\,\}$.  Let ${\bf m}\in C$ be given.  By the signature rule of tensor product of crystals with respect to $U_q(\mf{sl}_2)$ (cf. \cite[Remark 2.1.2]{KN}), we see that $|{\bf m}^{(2)}|-|{\bf m}^{(1)}|=n-k$. Let us identify each column ${\bf m}^{(b)}$ ($b=1,2$) with a semistandard tableau $T^{(b)}$ of  single-column  $(1^{|{\bf m}^{(b)}|})$ with entries in $[\ov{n}]$ such that $a\in [\ov{n}]$ appears in $T^{(b)}$ if and only if $m_{ab}=1$. Let $T$ be a tableau of two-column shape $(2^{|{\bf m}^{(1)}|},1^{|{\bf m}^{(2)}|-|{\bf m}^{(1)}|})$, whose left (resp. right) column is $T^{(2)}$ (resp. $T^{(1)}$). Then $T$ is semistandard by \cite[Lemma 6.2]{K13-2}. Since the actions of $\te_i$ and $\tf_i$ for $i\in I_0$ on $T$ coincides with those on ${\bf m}$ (see \cite[Section 7.1]{K13-2}), the map ${\bf m}\mapsto T$ gives a $U_q(\mf{g}_0)$-crystal (or $U_q(C_n)$-crystal) isomorphism  from $C$ to the set of semistandard tableaux of shape $(2^s,1^{n-k})$ ($0\leq s\leq k$) with entries in $[\ov{n}]$, and its image is isomorphic to $B_0({\rm cl}(\varpi_k))$ by \cite[Theorem 7.1]{K13-2}. This implies that $C=C({\bf v}_{k,0})$. Hence, we have $
B(W(\varpi_k))\cong\{\,{\bf m}\,|\,\sigma({\bf m})=(n-k,0)\,\}$.  

(2) We have $B(W(\varpi_k))\cong \bigsqcup_{l=0}^kC_0({\bf v}_{l,n-k})$ as a $U_q( C_n)$-crystal (see \eqref{classical decomposition-1}). Since $\td{F}^{n-k}{\bf v}_{l,n-k} ={\bf v}_{l,0}$ and $\td{F}$ commutes with $\te_i$ and $\tf_i$ for $i\in I_0$,  $\td{F}^{n-k} : C_0({\bf v}_{l,n-k}) \longrightarrow C_0({\bf v}_{l,0})$ is  a $U_q(C_n)$-crystal isomorphism. On the other hand, $C_0({\bf v}_{l,0})=\{\,{\bf m}\,|\,\sigma({\bf m})=(n-l,0)\,\}$ by (1), which implies that  $C_0({\bf v}_{l,n-k})=\{\,{\bf m}\,|\,\sigma({\bf m})=(k-l,n-k)\,\}$. Hence, we have 
$B(W(\varpi_k))\cong\{\,{\bf m}\,|\,\sigma({\bf m})=(k-l,n-k) \text{ for $0\leq l\leq k$}\,\}$.
 
(3) We have $B(W(\varpi_k))\cong C({\bf v}_{k,0})= C_0({\bf v}_{k,0})$. First, we have by the same argument as in Proposition \ref{decomposition of M} Case 2 that $|{\bf m}|\geq n-k$ for ${\bf m}\in C({\bf v}_{k,0})$.
Next, given ${\bf m}\in C({\bf v}_{k,0})$, we can check by \eqref{Kashiwara op for (1)} and induction on $|{\bf m}|$  that $\sigma({\bf m})=(n-k,t)$ for some $0\leq t\leq k$, which also implies that $|{\bf m}^{(2)}|=n-k+s$, $|{\bf m}^{(1)}|=s+t$ for some $s\geq 0$. 

Let $C=\{\,{\bf m}\,|\,\sigma({\bf m})=(n-k,t) \text{ for some $0\leq t\leq k$}\,\}$. Let ${\bf m}\in C$ be given with   $\sigma({\bf m})=(n-k,t)$ and $|{\bf m}^{(2)}|=n-k+s$, $|{\bf m}^{(1)}|=s+t$ for some $s,t\geq 0$.
Let $T^{(b)}$ be the semistandard tableau of single-column $(1^{|{\bf m}^{(b)}|})$ corresponding to ${\bf m}^{(b)}$ for $b=1,2$. Let $T$ be a tableau of skew shape $(2^{s+t},1^{n-k})/(1^{t})$, whose left (resp. right) column is $T^{(2)}$ (resp. $T^{(1)}$). Then the map ${\bf m}\mapsto T$ gives a $U_q(\mf{g}_0)$-crystal (or $U_q(B_n)$-crystal) isomorphism  from $C$ to a set of semistandard tableaux of shape $(2^{s+t},1^{n-k})/(1^t)$ ($s,t\geq 0$, $0\leq s\leq k$) with entries in $[\ov{n}]$ by \cite[Lemma 6.2]{K13-2}, and its image is isomorphic to $B_0({\rm cl}(\varpi_k))$ by \cite[Theorem 7.1]{K13-2}. This implies that $C=C({\bf v}_{k,0})\cong B(W(\varpi_k))$.  

(4) If $k=n$, then $B(W(\varpi_n))\cong \bigsqcup_{l=0}^{[\frac{n}{2}]} C_0({\bf v}_{n-2l ,0})$  as a  $U_q(C_n)$-crystal (see \eqref{classical decomposition-2}). Since $C_0({\bf v}_{n-2l,0})=\{\,{\bf m}\,|\,\sigma({\bf m})=(2l,0)\,\}$ by (1), we have 
$B(W(\varpi_n))\cong\{\,{\bf m}\,|\,\sigma({\bf m})=(2l,0) \text{ for $0\leq l\leq [n/2]$}\,\}$. 
If $k\neq n$, then we have
$
B(W(\varpi_k))\cong\{\,{\bf m}+\varsigma({\bf m})\,|\,\sigma({\bf m})=(2l,n-k) \text{ for $0\leq l\leq [n/2]$}\,\}$ by Corollary \ref{crystal realization} (see \eqref{2 fold decomposition of C}).  
\qed
\vskip 3mm

\begin{rem}{\rm
By \cite[Section 7.2]{K13-2}, we also obtain an explicit bijection from a classical crystal $C_0({\bf v}_{k,0})$ of type $B_n$ or $C_n$ to that of Kashiwara-Nakashima tableaux of (non-spinor) single column with length $k$. 
}
\end{rem}

%\appendix{}
%\section{Crystal graphs of $W(\varpi_k)$}

\end{document}